\documentclass[reqno,11pt]{amsart}
\usepackage{amssymb}
\usepackage{verbatim}
\newif\ifpdf
\ifpdf
  \usepackage[pdftex]{graphicx}
  \usepackage[pdftex]{hyperref}
\else
  \usepackage{graphicx}
\fi
\numberwithin{equation}{section}       
\theoremstyle{plain}
\newtheorem{Thm}{Theorem}[section]
\newtheorem{Prop}[Thm]{Proposition}
\newtheorem{Lemma}[Thm]{Lemma}
\newtheorem{Cor}[Thm]{Corollary}
\newtheorem{Prop-def}[Thm]{Proposition-Definition}

\newtheorem*{thmA}{Theorem A} 
\newtheorem*{thmB}{Theorem B} 

\theoremstyle{definition}

\newtheorem{Def}[Thm]{Definition}
\newtheorem{Remark}[Thm]{Remark}

\newtheorem*{setting}{Setting}

\newcommand{\C}{{\mathbf{C}}}

\newcommand{\N}{{\mathbf{N}}}
\newcommand{\Q}{{\mathbf{Q}}}
\newcommand{\R}{{\mathbf{R}}}
\newcommand{\Z}{{\mathbf{Z}}}

\newcommand{\fa}{{\mathfrak{a}}}
\newcommand{\fb}{{\mathfrak{b}}}
\newcommand{\fm}{{\mathfrak{m}}}

\newcommand{\cB}{{\mathcal{B}}}

\newcommand{\cH}{{\mathcal{H}}}
\newcommand{\cI}{{\mathcal{I}}}
\newcommand{\cJ}{{\mathcal{J}}}

\newcommand{\cO}{{\mathcal{O}}}

\newcommand{\cV}{{\mathcal{V}}}
\newcommand{\hcV}{\hat{\mathcal{V}}}

\newcommand{\elin}{{\mathcal L}^{\infty}}
\newcommand{\elde}{{\mathcal L}^2}

\newcommand{\fB}{{\mathfrak{B}}}

\newcommand{\fX}{{\mathfrak{X}}}

\newcommand{\hu}{{\hat{u}}}

\newcommand{\hvarphi}{{\hat{\varphi}}}

\renewcommand{\=}{:=}
\renewcommand{\a}{\alpha}

\newcommand{\e}{\varepsilon}

\newcommand{\Reg}{\mathcal{R}}

\newcommand{\eg}{e.g.\ }
\newcommand{\ie}{i.e.\ }
\newcommand{\cf}{cf.\ }

\newcommand{\cVqm}{{\mathcal{V}_{\mathrm{qm}}}}
\newcommand{\hcVqm}{{\hat{\mathcal{V}}_{\mathrm{qm}}}}
\newcommand{\hcVdiv}{{\hat{\mathcal{V}}_{\mathrm{div}}}}

\newcommand{\cVm}{{\mathcal{V}_{\mathrm{m}}}}
\newcommand{\cVdiv}{{\mathcal{V}_{\mathrm{div}}}}

\newcommand{\nuK}{\nu^\mathrm{K}}
\newcommand{\nuL}{\nu^\mathrm{L}}

\renewcommand{\div}{\operatorname{Div}}
\newcommand{\Div}{\operatorname{Div}}
\newcommand{\Cdiv}{\operatorname{C-Div}}
\newcommand{\ord}{\operatorname{ord}}

\newcommand{\Spec}{\operatorname{Spec}}
\newcommand{\Aff}{\operatorname{Aff}}
\newcommand{\PA}{\operatorname{PA}}
\newcommand{\MA}{\operatorname{MA}}

\newcommand{\trdeg}{\operatorname{tr.deg}}

\begin{document}
%
%

\setcounter{tocdepth}{1}

\title{Valuations and plurisubharmonic singularities.}
\date{\today}
\author{S{\'e}bastien Boucksom, Charles Favre \and Mattias Jonsson}

\address{CNRS-Universit{\'e} Paris 7\\
  Institut de Math{\'e}matiques\\
  F-75251 Paris Cedex 05\\
  France}
\email{boucksom@math.jussieu.fr}

\address{CNRS-Universit{\'e} Paris 7\\
         Institut de Math{\'e}matiques\\
         F-75251 Paris Cedex 05\\
         France\\ and
         Unidade Mista CNRS-IMPA\\
         Dona Castorina 110\\
          Rio de Janeiro\\
          22460-320\\
         Brazil}
\email{favre@math.jussieu.fr}

\address{Dept of Mathematics\\
  KTH\\
  SE-100 44 Stockholm\\
  Sweden\\
  and\\
  Dept of Mathematics\\
  University of Michigan\\
  Ann Arbor, MI 48109-1109\\
  USA}
\email{mattiasj@umich.edu}

\thanks{First author
  supported by the Japanese Society for the Promotion of Science.
  Third author supported by the NSF, the Swedish Research Council and 
  the Gustafsson Foundation.}

\dedicatory{Dedicated to Heisuke Hironaka on the occasion
of his seventy-seventh birthday}

\begin{abstract}
  We extend to higher dimensions some of the valuative analysis of
  singularities of plurisubharmonic (psh) functions developed by 
  the first two authors.
  Following Kontsevich and Soibelman we
  describe the geometry of the space $\cV$ of all 
  normalized valuations on $\C[x_1,\dots,x_n]$ centered at the 
  origin. It is a union of
  simplices naturally endowed with an affine structure. Using relative
  positivity properties of divisors living on modifications of
  $\C^n$ above the origin, we define formal psh functions on $\cV$,
  designed to be analogues 
  of 
  the usual psh functions. For
  bounded formal psh functions on $\cV$, we define a mixed
  Monge-Amp{\`e}re
  operator which reflects the intersection theory of divisors 
  above the origin of $\C^n$. 
  This operator associates to any $(n-1)$-tuple of formal psh 
  functions a positive measure of finite mass on $\cV$.
  Next, we show that the collection of Lelong numbers of a 
  given germ $u$ of a psh function at all infinitely near points 
  induces a formal psh function
  $\hat u$ on $\cV$. When $\varphi$ is a psh H{\"o}lder weight in
  the sense of Demailly, the generalized Lelong number
  $\nu_\varphi(u)$ equals the integral of $\hu$ against the
  Monge-Amp{\`e}re measure of $\hvarphi$.
  In particular, any generalized Lelong number is an 
  average of valuations.
  We also show how to compute the multiplier ideal of $u$ and the 
  relative type of $u$ with respect to $\varphi$ in the sense of
  Rashkovskii, in terms of $\hat u$ and $\hat\varphi$.
\end{abstract}


\maketitle


\tableofcontents

\newpage
%
%
%
%
\section*{Introduction}
This paper aims at laying down the foundations for a valuative study
of singularities of plurisubharmonic (psh) functions in any
dimension. Such an analysis was carried out by the last two authors in
dimension two in a series of
papers~\cite{valtree,FJ-sing,FJ-multsing}. But that work used in an
essential way the assumption on the dimension, and did not immediately
extend to a more general setting. Here we present a more geometric
approach, using as a key tool Hironaka's theorem on resolution of
singularities

\smallskip
Fix $n\ge1$ and let $u:(\C^n,0)\to \R\cup\{ - \infty\}$ be a 
\emph{psh germ}, 
that is, a psh function defined in a neighborhood of the origin. 
We say that $u$ has a \emph{singularity} at $0$ when $u(0) = -\infty$. 
Various quantities have been introduced to measure the singularity 
of $u$ at $0$, the most classical certainly being the so-called 
\emph{Lelong number}
$\nuL(u,0)=\max\{c >0, \, u(z) \le c \,\log | z | + O(1) \}$. 
Although of fundamental importance, the Lelong number gives 
too rough an insight
into the singularity of $u$ in many situations.  To remedy this
problem, we propose to look at the family of Lelong numbers
$\nuL(u\circ \pi,p)$ at all points $p\in \pi^{-1}(0)$ lying in some
birational model $\pi: X_\pi \to X\=\C^n$.  The general aim of this
article is to show that given this collection of Lelong numbers, one
can recover essentially all other known measurements of the of
singularity of $u$, including multiplier ideals, 
generalized Lelong numbers and relative types.
\begin{thmA}
  Let $u,v: (\C^n,0) \to \R\cup\{-\infty\}$ be two psh germs.
  Consider the following four statements.
  \begin{enumerate}
  \item
    For all proper modifications $\pi: X_\pi \to \C^n$ above 
    $0$ and all points $p\in\pi^{-1}(0)$, 
    we have $\nuL(u\circ \pi,p) = \nuL(v\circ \pi,p)$.
  \item
    For all $t>0$ we have 
    $\cJ_+(tu)=\cJ_+(tv)$.
  \item 
    For any tame maximal psh weight $\varphi$, the relative types
    $\sigma(u,\varphi) =\sup\{ c>0,\, u \le c\, \varphi + O(1) \}$ 
    and $\sigma(v,\varphi)$
    are equal.
  \item
    For any tame psh weight $\varphi$,  
    the generalized Lelong numbers 
    $\nu_\varphi(u) = dd^c u \wedge (dd^c\varphi)^{n-1} \{ 0 \}$
    and $\nu_\varphi(v)$ are equal.
  \end{enumerate}
  Then (1), (2) and (3) are equivalent, and imply (4).
\end{thmA}
Let us explain the notation used in the statement of the theorem. For
any psh germ $u$, we denote by $\cJ(u)$ its \emph{multiplier ideal}, \ie
the ideal of germs of holomorphic functions $f$ at $0$ 
such that $fe^{-u}\in L^2_\mathrm{loc}$. 
In (2) $\cJ_+(u)$ refers to the stationary value of
$\cJ((1+\e)u)$ as $\e>0$ tends to $0$.
    
A \emph{psh weight} is a psh germ $\varphi$ with an isolated
singularity at the origin such that $e^\varphi$ is continuous. 
A psh weight is said to be \emph{tame} if
it can be strongly approximated by its multiplier ideals, \cf
Definition~\ref{def:reg}. In particular psh weights $\varphi$ such
that $e^\varphi$ is H{\"o}lder continuous are tame.  Finally, a weight
is said to be \emph{maximal} if it is locally maximal outside the
origin, \ie its Monge-Amp{\`e}re measure is a Dirac mass at $0$.  The
notion of generalized Lelong number is due to Demailly~\cite{dem},
whereas the relative type was introduced more recently by
Rashkovskii in~\cite{rash}.  The theorem above shows that these
notions can be read off directly from the Lelong numbers of $u$ in all
birational models.  We expect (4) to imply in turn the other
statements and that $\cJ_+$ can be replaced by $\cJ$ itself in~(2) but
further analysis is required to prove this.

\[\diamond\]

At the heart of our method lies the idea of organizing the
data $\{\nuL(u,p)\}$ in an efficient way, and turning it into a function
on a suitable space of valuations. We shall prove that both this
function and this space exhibit striking properties which relate them
in a natural way to convex geometry. This feature will eventually help
us to control the behavior of the Lelong numbers $\nuL(u\circ \pi,p)$
when moving the point $p$ in the tower of all blow-ups $\pi$.

To motivate the introduction of valuations, recall that
for any holomorphic function $f$, the integer $\nuL(\log|f|)$ equals the
multiplicity of $f$ at $0$. More generally, for any modification $\pi:
X_\pi \to \C^n$, and any point $p\in\pi^{-1}(0)$, the function $f
\mapsto\nuL(\log|f\circ\pi|,p)$ is equal to the divisorial valuation
$f\mapsto\nu_p(f)\= \ord_p (f \circ \pi)$, 
defined on the ring $R\= \C[x_1,\dots,x_n]$.  If $\cVdiv$ denotes the space of all (suitable normalized) divisorial valuations
$\cVdiv$ centered at $0$, we see that the collection $\{ \nuL(u\circ
\pi,p)\}$ induces a real-valued function on $\cVdiv$. As we
shall soon explain, this function, which will be denoted by $-\hu$, 
is of a very special type.

Unfortunately, the natural topology of pointwise convergence turns
$\cVdiv$ into a totally discontinuous space, and it is difficult to
directly study functions on it. One is led to complete this space in
some way. As it turns out, there are two reasonable
choices that one can make, which we refer to as Zariski's and
Berkovich's point of view.  As we shall see, both approaches shed
their own light onto the problem, and benefit one from another.

In Zariski's perspective~\cite{Z}, one considers $\fX \=
\varprojlim_\pi X_\pi$ the projective limit over all proper
modifications of $\pi: X_\pi \to X \= \C^n$ above $0$.
This space contains
all rank $1$ valuation rings of the function field $\C(X)$ whose
maximal ideal contain the maximal ideal $\fm$ of $X$ at $0$. We shall 
refer to it as the \emph{Riemann-Zariski space} of $X$ at $0$. It clearly contains $\cVdiv$, and admits a natural topology 
for which it
is quasi-compact. In $\fX$, divisors are in one-to-one correspondence
with functions on the subspace $\cVdiv \subset \fX$, hence any psh
germ defines a divisor $Z(u)$ on $\fX$.

The Berkovich point of view refers to the theory of analytic spaces
over non-archimedean fields as developed in~\cite{ber}.  In our
context, we view $\C$ as a field with a trivial valuation, and define $\cV$ to be the set of all valuations $\nu:R\to
[0,+\infty]$ that extend this trivial valuation on $\C$, 
are $>0$ on the maximal
ideal $\fm=(x_1,\dots,x_n)\subset R$ and normalized by
$\nu(\fm):=\min\nu(x_i)=1$.  From this perspective, $\cVdiv$ becomes
a dense subset of $\cV$, which is compact for the topology of
pointwise convergence. Of course, given a psh germ $u$, it is
unclear whether the function $\hat u$, which is a priori defined only on $\cVdiv$ extends in a natural way to $\cV$.  
But we may remark that in
the case $u=\log|f|$ for some $f\in R$, the function $\hat u(\nu):=-\nu(f)$
is a natural choice, as it yields a continuous function on $\cV$.  We
shall prove that for any psh germ $u$, the function $\hat u$ does
indeed extend to $\cV$ in a unique way with an appropriate 
continuity property.

\[\diamond\] 

To go further, it is now necessary to take a closer look at the
space of divisors on $\fX$, and at the geometry of $\cV$. Let us begin
with the latter, following~\cite{KKMS} and~\cite{KoSo}.

If $\pi: X_\pi \to \C^n$ is a modification above 
the origin with simple
normal crossing exceptional divisor, then one can construct the 
\emph{dual} (simplicial) \emph{complex} 
$\Delta(\pi)$ of $\pi$ which encodes the incidence
relations between the irreducible components of $\pi^{-1}(0)$. 
Its underlying topological space is a
union of simplices of real dimension $\le n-1$, 
and is homotopic to a point. 
What makes this
construction useful is the fact that vertices of $\Delta(\pi)$ are
irreducible components of $\pi^{-1}(0)$, and can thus be viewed as
divisorial valuations included in $\cV$.  
Moreover, the faces of $\Delta(\pi)$ can be naturally
identified with valuations that are monomial in suitable local
coordinates in $X_\pi$. 
Thus we can realize the underlying topological space of $\Delta(\pi)$
as a subset $|\Delta(\pi)|$ of $\cV$. 
When $\pi'$ dominates $\pi$, \ie when there exists a
modification $\mu: X_{\pi'} \to X_\pi$ such that $\pi'=\pi\circ\mu$,
we have a natural inclusion $|\Delta(\pi)|\subset|\Delta(\pi')|$.
Letting $\pi$ tend to infinity in the net of all
modifications yields a dense inclusion 
from $\cup_\pi|\Delta(\pi)|$ into $\cV$.

Being identified with sets of monomial valuations, faces of dual
complexes are naturally parameterized by the weights on the
corresponding coordinates in which the valuations are monomial 
and thus inherit an \emph{affine structure}.  
A \emph{simplex} in $\cV$ is a face of
some $|\Delta(\pi)|$ endowed with this affine structure.  
The affine structure on $\cV$ is of fundamental importance.

Returning to a psh germ $u$, it is not difficult to interpret the
restriction of $\hat u$ to a fixed simplex in terms of a weighted version of Lelong numbers called 
\emph{Kiselman numbers}~\cite{kis}. Using this fact,
we conclude that the restriction of $\hat u$ to any 
simplex is convex. This 
is already a quite remarkable fact, but it is
also of very local nature as simplices are relatively small subsets of
$\cV$.  One may thus wonder if the function $\hat u$ does not exhibit
more global convexity properties. Unfortunately it is not immediately
clear how all the simplices patch together to build up $\cV$,\footnote{Except
when $n=2$ where simplices are segments, and $\cV$ is a real tree, see~\cite{valtree}.}
so it seems hard to define any notion of convexity in $\cV$
starting from the simplices.  To overcome this difficulty, we need to
turn our attention to Zariski's approach, describe the space of
divisors on $\fX$ and study positivity properties on this space in the
spirit of what was done in~\cite{BFJ}.

\smallskip
A \emph{Weil divisor} $Z$ on $\fX$ is defined as a collection $\{Z_\pi\}$
of divisors in all modifications $X_\pi \to X$ above $0$,
compatible under push-forward. \emph{Cartier divisors} are those 
that are determined at a finite level $X_\pi$, that is, 
$Z_{\pi'} = \mu^*Z_\pi$ 
for any $\mu: X_{\pi'} \to X_\pi$. Clearly, a Weil divisor $Z$
is the same as a function $g_Z$ on $\cVdiv$, and it is not difficult
to see that when $Z$ is Cartier, $g_Z$ extends as a continuous
function on $\cV$.  The space $\Cdiv(\fX)$ of Cartier divisors sits as
a dense subspace inside the space $\div(\fX)$ of Weil divisors, and
both spaces are infinite dimensional.

Now any psh germ $u$ induces a Weil divisor $Z(u)$ given in $X_\pi$ by
the Siu decomposition $dd^c (u\circ \pi) = -[Z(u)_\pi] + T'_\pi$ with
$T'_\pi$ a positive closed $(1,1)$ current not charging $\pi^{-1}(0)$.

A closer look reveals that the divisor $Z(u)$ possesses a
quite strong numerical positivity property: it lies in the closure of
the cone of nef Cartier divisors, \ie those Cartier divisors on $\fX$
that are determined by a nef divisor on some $X_\pi$. A Weil divisor
in this closure will be said to be \emph{nef}. We shall prove that 
for any
nef Weil divisor $Z$ the corresponding function $g_Z$ on $\cVdiv$
extends to $\cV$ in a natural way, see Section~\ref{sec:approx}.

Using this notion of nef Weil divisor, we finally arrive at the global
convexity notion on $\cV$ that we were seeking. 
Namely, we define a \emph{formal psh function} to be
an upper semicontinuous function $g:\cV\to[-\infty,0]$ that is
continuous on each simplex and whose restriction to $\cVdiv$ defines
a nef Weil divisor.
A formal psh function is convex on each simplex in $\cV$, 
but the global convexity property recasts informations on the differential of the function when
passing from one simplex to another. 

As we shall see, formal psh functions on $\cV$
have many properties in common with psh germs in $\C^n$.
By construction, any function $\hat u$ associated to a psh germ 
is formal psh on $\cV$. We call it the \emph{valuative transform}
of $u$.

In the philosophy underlying~\cite{BGS},
one might want to think of a nef Weil 
divisor $Z$ as the formal closed, positive $(1,1)$-current 
$i\partial\overline{\partial}g_Z$ associated to the formal psh function $g_Z$. 

\[\diamond\] 

We have thus associated to a psh germ $u$ on $(\C^n,0)$ a
nef Weil divisor $Z(u)$ on $\fX$ and a formal psh function 
$\hat u$ on $\cV$.
They encode the same information but from
two different perspectives. 
Let us now indicate how one can use these notions to actually prove
Theorem~A.

Condition~(1) says that two psh germs $u$ and $v$ have the same
associated valuative transforms: $\hat u=\hat v$. 
Thus the equivalence between~(1)
and~(2) follows from the understanding of the connection between 
$\hat u$
and the multiplier ideals of $u$. We put this problem into a larger
perspective and define for any formal psh function $g$ on $\cV$ a
\emph{valuative multiplier ideal} denoted by $\elde(g)$, which is an ideal of
the ring $\hat R=\C[[x_1,\dots,x_n]]$ of formal power series. Although
we do not know\footnote{If true, this would prove the Demailly-Koll{\'a}r Openness Conjecture,
   see below.}  if the valuative multiplier ideal $\elde(\hu)$
coincides with the usual multiplier ideal. $\cJ(u)$ in the ring $\cO$
of holomorphic germs at $0$, we can prove that 
$\elde(\hat u)\cap\cO=\cJ_+(u)$, and this
shows that~(1) implies~(2). Conversely, Demailly's approximation
results show that the singularities of a psh germ $u$ can be
approximated by its multipliers ideals $\cJ(ku)$ as $k\to\infty$ in a
very precise way, and we show that this is enough to recover 
$\hat u$, extending the analysis of~\cite{FJ-multsing}.
Thus~(2) implies~(1).

That~(1) implies~(3) is a consequence of the following result: if
$\varphi$ is a tame psh weight, then a psh germ $u$ satisfies
$u\leq\varphi+O(1)$ iff $\hat u\leq\hat\varphi$.  This is proved by
another application of Demailly's results. Conversely, we show 
that~(3) implies~(1) by showing that quasi-monomial valuations can be
represented as relative types with respect to tame psh weights.

To prove that~(1) implies~(4) we compute the mass of the
Monge-Amp{\`e}re measure $dd^cu \wedge (dd^c\varphi)^{n-1}$ 
at the origin in terms of the formal psh functions $\hat u$ 
and $\hat\varphi$ (with
$\varphi$ a psh weight). 
This is done by introducing an intersection
theory on nef divisors on $\fX$ and transferring its outcome to $\cV$.
More precisely, we proceed as follows.

The intersection of $n$ divisors in a fixed model $X_\pi$ leads to the
definition of the intersection product $\langle Z_1,\dots,Z_n \rangle$
of $n$ Cartier divisors. Using a monotonicity property of this
product with respect to nef divisors, we extend the definition to
arbitrary nef Weil divisors by approximation, and give a meaning to
the intersection $\langle Z, Z_2,\dots, Z_n\rangle$ of a Cartier
divisor $Z$ against nef Weil divisors. We then show that the linear
map $g_Z \mapsto \langle Z, Z_2,\dots, Z_n\rangle\in \R$ extends to a
positive linear continuous functional on all continuous functions on
$\cV$. This functional is thus attached to a positive Radon measure,
which we call the 
\emph{Monge-Amp{\`e}re measure} of the formal psh functions $g_i$
associated to the nef Weil divisors $Z_i$. We denote it by
$\MA(g_1,\dots, g_{n-1})$.  As a measure on the $(n-1)$-dimensional
space $\cV$ it formally possesses many features of its complex
analytic analog $dd^c u_1\wedge\dots\wedge dd^cu_n$ on $\C^n$. The
statement (1) $\Rightarrow$ (4) reflects a deeper connection between
these two theories.
\begin{thmB}
  For any tame psh weight $\varphi$ and any psh germ $u$, the
  formal psh function $\hat u$ is integrable with respect to 
  the positive Radon measure $\MA(\hat\varphi)$, and
  \begin{equation*}
    \nu_\varphi(u)
    \=(dd^c u \wedge (dd^c\varphi)^{n-1})\{ 0 \} 
    =\int_{\cV} -\hat u\, \MA(\hat\varphi).
  \end{equation*}
  In particular, the generalized Lelong number $\nu_\varphi$
  is an average of valuations.
\end{thmB}
Here $\MA(\hvarphi)$ is a shorthand for 
$\MA(\hvarphi,\dots,\hvarphi)$.
Technically the proof relies on yet 
another application of Demailly's approximation technique.  

\[\diamond\] 

The present work can be seen as an instance of developing a 
(pluri)potential theory in a non-Archimedean setting. For
related works, see \eg~\cite{autissier,BR1,BGS,bost,Ch,CLR,gubler,kani,maillot,thuillier1,thuillier2,zhang}. 
There are several questions left open by the present work, 
notably the remaining implication~(4)$\implies$(1) in Theorem~A.
One way to prove this implication would be to find 
fundamental solutions to the formal Monge-Amp{\`e}re operator.
Another interesting question is the Openness Conjecture by
Demailly and Koll{\'a}r~\cite[Remark~5.3]{DK}: 
is $\cJ_+(u)=\cJ(u)$ for any psh germ $u$? 
The technique used in~\cite{FJ-multsing} to prove this conjecture
in dimension two relies on a more detailed understanding
on the behavior of formal psh functions. 
It would also be of interest to generalize
the dynamical results from~\cite{FJ-eigenval} to
higher dimensions.
We plan to come back to these issues in future work.

\[\diamond\] 

\noindent{\bf Organization of the paper.} 
 Section~\ref{sec:geom}
contains basic definitions concerning the space of divisors in $\fX$,
and a thorough discussion of the geometry of $\cV$, mostly borrowed
from~\cite{KKMS} and~\cite{KoSo} with a few adaptations to our 
local context. In
Section~\ref{sec:pos}, we undertake the study of positivity properties
of divisors in $\fX$ and transport them to $\cV$ with its affine
structure.  Section~\ref{sec:mult} contains a general discussion of
valuative multiplier ideals.  These are used to approximate 
formal psh functions by ideals, a powerful technique that 
plays a fundamental
role in the definition of the Monge-Amp{\`e}re operator, see
Remark~\ref{rem:essential}.  The latter operator is introduced and
discussed in Section~\ref{sec:inter}. Finally we prove Theorems~A
and~B in Section~\ref{sec:psh}.

\begin{setting}
In Sections~\ref{sec:geom} to~\ref{sec:inter}, we let $k$ be an
algebraically closed field of characteristic $0$ and choose $0\in X$
as a smooth closed point in an $n$-dimensional affine $k$-variety,
with coordinate ring $R$ and maximal ideal $\fm$ at $0$. The formal
completion of $X$ at $0$ will be denoted by $\hat X$. It is the
infinitesimal neighborhood of $0$ in $X$, and is described by its
coordinate ring $\hat R$, the completion of $R$ in the $\fm$-adic
topology. Since $0$ is a smooth point, $\hat R$ is isomorphic to
$k[[x_1,\dots,x_n]]$ for any choice of a coordinate system $x_1,\dots,x_n$ at $0$.

In Section~\ref{sec:psh} where we speak about psh functions, we shall
return to the com\-plex setting $k=\C$, with $0$ being the origin in
$X=\C^n$.
\end{setting}

\medskip

The following table might help the reader in translating notions from 
the geometric (or Zariski) point of view to the functional (or Berkovich) one.

\begin{center}
\begin{tabular}{|l|l|}
\cline{1-2}
Berkovich  point of view (functional) & Zariski point of view (geometric)
\\
\cline{1-2}
Valuative space $\cV$ \S\ref{sec:val}
&
Riemann-Zariski space $\fX$ \S\ref{sec:RZ}
\\
Functions on $\cV_\mathrm{div}$ 
&
Weil divisors \S\ref{sec:divisors}
\\
Formal psh functions \S\ref{sec:pos1}
&
Nef Weil divisors \S\ref{sec:convx}
\\
Monge-Amp{\`e}re operator \S\ref{sec:inter2}
&
Intersection theory of nef Weil divisors \S\ref{sec:inter1}
\\
\cline{1-2}
\end{tabular}
\end{center}

%
%
%
%
%
\section{Geometry of valuation spaces}\label{sec:geom}
%
%
\subsection{The Riemann-Zariski space}\label{sec:RZ}
This space was introduced by Zariski~\cite{Z} as a tool for
desingularizing algebraic varieties, see~\cite{ZS,vaquie,Co} for a
more detailed account on its definition and its history. The ringed
space structure on the Riemann-Zariski space was introduced by
Hironaka in~\cite{Hir}. 

Denote by $\fB$ the set of all projective
modifications $\pi: X_\pi \to X$ above the origin 
(\ie $\pi$ is an isomorphism above $X\setminus\{0\}$)
such that $X_\pi$ is smooth. If $\pi, \pi'
\in \fB$, we say that $\pi'$ \emph{dominates} $\pi$ and write $\pi'\ge
\pi$, if there exists a morphism (necessarily unique and birational)
$\mu:X_{\pi'}\to X_\pi$ such that $\pi'=\pi\circ\mu$. This endows
$\fB$ with a partial order relation. By Hironaka's desingularization
theorem, this ordered set forms a directed family.  The
\emph{Riemann-Zariski space} of $X$ at $0$ is the projective limit of
locally ringed topological spaces
\[\fX:=\varprojlim_\pi X_\pi,\]
where $X_\pi$ is viewed as a scheme. It can be shown that the locally ringed 
space $\fX$ is not isomorphic to a scheme, even of infinite type. 
\begin{Remark}\label{rem:RZ}
A comment on our terminology is in order.
In dimension $2$, the space $\fX$ coincides with the set $\hat{\fX}$ of all valuation rings
containing $\fm$, but this fact is no longer true as soon as $n\ge 3$.
There are however natural maps  $\imath :\fX \hookrightarrow \hat{\fX}$, and $p: \hat{\fX} \to \fX$ such that $\imath$ is injective, $p$ is surjective and $p\circ \imath$ is the identity on $\fX$.
The map $p$ is defined by sending a valuation to the family of its centers in all models $X_\pi$ with $\pi \in \fB$. Note that Zariski's theorem says that a valuation is determined by the family of its centers in \emph{all} birational models, not only those which are isomorphic outside the origin.
For a point $x = \{ x_\pi\} \in \fX$, let $\cI_{x}$ be the prime ideal of those functions $f\in R$
for which $x_\pi$ is included in the strict transform of $f=0$ in $X_\pi$ for all $\pi$. 
When $\cI_{x}= (0)$ the union of all local rings $\cO_{x_\pi}$ is a valuation ring which we define to be  $\imath(x)$. Otherwise,
define $C_{x,1} = \mathrm{Spec}\, \cO/\cI_x$.
The collection $\{x_\pi\}$ now defines a point in the projective limit of the strict transforms of $C_{x,1}$
by all $\pi$. We can thus inductively define a nested sequence of irreducible subschemes
$C_{x,1} \supsetneq C_{x,2} \supsetneq \cdots C_{x,k}$. On $\cO_{C_{x,k}}$ the point $x$ defines a valuation ring $\nu$.
Set $\imath(x)$ to be the composed valuations of the divisorial valuations $\ord_{C_{x,j}}$ for $j=1, ..., k-1$ together with $\nu$.
From this description of $\imath$ and $p$, one sees that $\fX$ contains all rank $1$ valuations centered at $0$.
\end{Remark}

%
%
\subsection{Divisors on $\fX$}\label{sec:divisors}
On the projective limit space $\fX$, one can define in a natural way
Weil and Cartier divisors as limits of Weil and Cartier divisors on
the $X_\pi$'s. This is basically Shokurov's notion of birational
divisors (or b-divisors) on $X$, \cf~\cite{shokurov}. The subsequent
discussion is adapted from~\cite{BFJ}.

\smallskip
For $\pi\in \fB$, we write $\div(\pi)$ for the set of divisors on
$X_\pi$ with real coefficients and exceptional for $\pi$. It is a
finite dimensional real vector space endowed with its canonical
topology. Note that an exceptional divisor is uniquely determined by
its numerical class, so that we obtain a natural isomorphism between
$\div(\pi)$ and the relative N{\'e}ron-Severi space $NS(X_\pi/X)_\R$ of
$X_\pi$ (recall that $X_\pi$ is assumed to be smooth).  
A birational morphism $\mu : X_{\pi'}\to X_\pi$ 
with $\pi', \pi \in \fB$ induces natural
linear maps 
$\mu^*: \div(\pi) \to \div(\pi')$, $\mu_* : \div(\pi') \to\div(\pi)$.
\begin{Def}{~}
\begin{itemize}
\item
The space of \emph{Weil divisors} on $\fX$ is the projective limit
\[\div(\fX) \=\varprojlim_\pi \div(\pi)\]
with respect to the push-forward arrows. It is endowed with its
projective limit topology, which will be called the \emph{weak
topology}.
\item
The space of \emph{Cartier divisors} on $\fX$ is the inductive limit
\[\Cdiv(\fX) :=\varinjlim_\pi \div(\pi)\]
with respect to the pull-back arrows.
It is endowed with its
inductive limit topology, which will be called the 
\emph{strong topology}.
\end{itemize}
\end{Def}
Note that these vector spaces are infinite dimensional since
we assume $n\ge2$. 
We will denote by $\Div(\fX)_\Z$ and $\Cdiv(\fX)_\Z$ the
spaces of Weil and Cartier divisors with $\Z$-coefficients.

A Weil divisor $Z\in \div(\fX)$ is by definition described by its
\emph{incarnations} $Z_\pi\in\div(\pi)$ on each smooth
birational model $X_\pi$ of $X$, compatible with each other by
push-forward. A sequence (or a net) $Z_i$ converges to $Z$ in
the weak topology iff $Z_{i,\pi} \to Z_\pi$ in
$\div(\pi)$ for each $\pi$.

On the other hand, the relation $\mu_*\mu^*\alpha=\alpha$
when $\mu$ is a birational morphism shows that there is an injection
\[\Cdiv(\fX)\hookrightarrow\div(\fX)\]
\ie a Cartier divisor is in particular a Weil
divisor.  Concretely, a Weil divisor $Z$ is Cartier iff there exists
$\pi$ such that its incarnations $Z_{\pi'}$ on higher blow-ups
$X_{\pi'}$ are obtained by pulling-back $Z_\pi$. We will call
such a $\pi$ a \emph{determination} of $Z$. There are also 
natural injective maps
\[
\div(\pi)\hookrightarrow \Cdiv(\fX)
\]
which extend a given divisor
$Z\in\div(\pi)$ to a Cartier divisor by pulling it back. This Cartier
divisor is by construction determined by $\pi$. In the sequel, we
shall always identify a divisor $Z \in \div(\pi)$ with its image in
$\Cdiv(\fX)$.

A function on $\Cdiv(\fX)$ is continuous in the strong topology iff
its restriction to each $\div(\pi)$ is continuous.
In particular, the natural injective map
$\Cdiv(\fX)\to\div(\fX)$ 
is continuous, and moreover has dense image.
Indeed, if $Z$ is a Weil divisor on $\fX$, then the incarnations 
$Z_\pi$ of $Z$ on $X_\pi$, viewed as Cartier (hence Weil) 
divisors on $\fX$, tautologically converge to $W$ as 
$\pi\to\infty$.

If $\fa\subset \hat R$ is an ideal, we let $Z(\fa)$ be the Weil
divisor on $\fX$ such that $\ord_E Z(\fa)=-\ord_E(\fa)$ for each
exceptional prime $E$. Note that $Z(\fa)\leq 0$. When $\fa=(f)$
is principal, we also write $Z(f)$. 
The divisor
$Z(\fa)$ on $\fX$ is Cartier iff $\fa$ is a primary ideal. By the
Nullstellensatz, a primary ideal $\fa$ of $\hat R$ is generated by
$\fa\cap R$, so that primary ideals in $R$ and $\hat R$ coincide. In
that case, $Z(\fa)$ is the Cartier divisor such that $-Z(\fa)_\pi$ is
the effective Cartier divisor defined by the principal ideal sheaf
$\pi^{-1}\fa$ whenever $\pi$ dominates the blow-up of~$\fa$. 
%
%
\subsection{Valuations}\label{sec:val}
By a \emph{valuation} $\nu$ on $X=\Spec(R)$, we mean a function
$\nu:R\to[0,+\infty]$ such that $\nu(f_1f_2)=\nu(f_1)+\nu(f_2)$,
$\nu(f_1+f_2) \ge \min \{ \nu(f_1) , \nu(f_2) \}$ for all $f_1, f_2
\in R$, and which induces the trivial valuation on $k\subset R$, \ie
$\nu(0)=+\infty$, $\nu|_{k^{*}}\equiv0$, and satisfying $\nu(\fm)
>0$.  This is the same as a non-Archimedean absolute value
$\exp(-\nu)$ on $R$ that restricts to the trivial one on $k$. We will
denote by $f\mapsto |f|(\nu):=\exp(-\nu(f))$, or simply
$f\mapsto|f|$, the corresponding
absolute value, following standard notation in non-Archimedean
geometry.

The ideal $\nu^{-1}(+\infty)$ is a prime ideal, cutting out a subvariety
$H=H(\nu,X)$ of $X$ called the \emph{home} of $\nu$ in $X$. Note that
$\nu$ induces a real-valued valuation in the sense of~\cite[\S8]{ZS}
on the function field $k(H)$ of its home, so we can attach
to $\nu|_H$ the usual invariants of valuation theory. In particular,
we can consider the value group
$\Gamma_\nu:=\nu(k(H)^{*})\subset\R$ and residue field
$\kappa_\nu:=\{\nu\geq 0\}/\{\nu>0\}$. Then the
Abhyankar inequality holds:
\begin{equation}\label{e:abh}
  \trdeg_k( \kappa_\nu)+\dim_\Q( \Gamma_\nu\otimes_\Z \Q) \leq\dim
  H(\nu,X).
\end{equation}

The \emph{center} of a valuation $\nu$ on $X$ is the subvariety of $H$
defined by the prime ideal $\{\nu>0\}\subset R$.

There exists a unique valuation whose home is reduced to $0$, namely
the trivial valuation $\nu_0$ sending all elements of $\fm$ to
$+\infty$.  Now pick any non-trivial valuation $\nu$.  For any
projective modification $\pi:X_\pi\to X$ of $X$ above $0$, $\nu$
induces a valuation on $X_\pi$, \ie a collection of compatible
valuations on the local rings of $X_\pi$. This is clear when the home
of $\nu$ is $X$ itself since $X$ and $X_\pi$ have the same function
fields. And the general case reduces to this one: the home of $\nu$ on
$X_\pi$ is then the strict transform of its home on $X$. In
particular, we can consider the center of $\nu$ on $X_\pi$: it is an
irreducible subvariety.

If $\fa$ is an ideal of $R$, we set
$\nu(\fa):=\min\{\nu(f),\,f\in\fa\}$, so that 
$\log|\fa|(\nu)=-\nu(\fa)$. 
A valuation $\nu$ on $R$ is centered at $0$ iff
$\nu(\fm)>0$, which holds iff the corresponding absolute value
$|\cdot|(\nu)$ continuously extends to an absolute value on 
$\hat R$ with its $\fm$-adic topology. 
Thus the set of valuations centered at $0$ on $X$ 
can be identified
with the set of continuous valuations on $\hat R$. Any non-trivial
valuation $\nu$ on $X$ centered at $0$ satisfies $0<\nu(\fm)<+\infty$,
so that we can normalize it as follows.
\begin{Def}
  We denote by $\hcV$ the set of all valuations $\nu$ on $X$ centered
  at $0$, and by $\cV$ the subspace of those valuations
  normalized by $\nu(\fm)=1$.  We endow $\hcV$ and $\cV$ with the
  topology of pointwise convergence on $R$. 
 \end{Def}
Then $\cV$ is compact, and $\hcV$ is endowed with a $\R_+^*$-action by scaling. The trivial valuation $\nu_0$ is the unique fixed point, and for each $\nu\in\hcV$ we have $t\nu\to\nu_0$ as $t\to+\infty$.  The space $\hcV\setminus \{\nu_0\}$ can be 
identified with the open cone over $\cV$. 
We will say that a function $g$ on $\hcV\to[-\infty,+\infty]$ is 
\emph{homogeneous}
 if it satisfies $g(t\nu)=tg(\nu)$ for each $t>0$ and each $\nu\in\hcV$.  
\begin{Remark}
  The space $\hcV$ we consider here is by definition the
  \emph{$k$-analytic space} in the sense of Berkovich associated by
  Thuillier in~\cite{thuillier2} to the formal completion $\hat X$ of the
  $k$-scheme $X$ at the origin. The singleton $\{\nu_0\}$ is the 
  $k$-analytic   space associated to the point space $\{0\}$, and the open cone
  $\hcV\setminus\{\nu_0\}$ over $\cV$ is the
  \emph{generic fiber} of $\hat X$.
\end{Remark}
%
%
\subsection{Divisorial and quasimonomial valuations}\label{sec:quasi}
A \emph{divisorial valuation} on $X$ is a valuation $\nu$
proportional to
the vanishing order $\ord_E$ along a prime divisor $E$ in some smooth
birational model $Y$ of $X$.  If $\nu$ is centered at $0$, then
$Y$ can be chosen as $X_\pi$ for some $\pi\in\fB$ and $E$ as an
\emph{exceptional prime} (divisor), that is, an irreducible component of
$\pi^{-1}(0)$.

A \emph{quasi-monomial} valuation $\nu$ on $X$ is a valuation that
becomes \emph{monomial} on some smooth birational model $Y$ of
$X$. This means the following: there exist a closed point $p$ in $Y$,
local coordinates $z=(z_1,\dots,z_n)$ at $p$ and real numbers
$w_1,\dots,w_n\ge0$ (not all zero) called weights such that if $f$ in
the local ring of $Y$ at $p$ is written as a power series $f=\sum_\a
a_\a z^\a$ in multi-index notation, then $\nu$ coincides with the
monomial valuation 
\[\nu_{z,w}(f):=\min\{w\cdot\a, a_\a\ne0\}.\]
In terms of the associated absolute value, a valuation is monomial in the
$z$-coordinates iff $|\sum_\a a_\a z^\a| =\max_\a |a_\a z^\a |$.  The
center of $\nu$ on $Y$ is then the intersection of all $\{z_i=0\}$ for
which $w_i>0$.  In particular, a divisorial valuation is
quasi-monomial (take $w_i=0$ for all but one $i$) and conversely $\nu$
as above is divisorial iff $w=(s_1,\dots,s_n)\in\R_+^n$ is
proportional to an element of $\Q_+^n$.  The corresponding prime
divisor is then the exceptional divisor of a suitable toric blow-up of
$Y$ in the coordinates $z_i$.

Quasi-monomial valuations can also be characterized as \emph{Abhyankar
valuations}, \ie real-valued valuations on $X$ for which equality
holds in the Abhyankar inequality~\eqref{e:abh},
see~\cite[Proposition~4.8]{ELS}.

If $\nu$ is a quasi-monomial valuation on $X$ centered at 0,
then the birational model $Y$ as above can be chosen to be $X_\pi$ for
some $\pi\in\fB$.  More precisely, let $\pi:X_\pi\to X$ be a
modification above 0 such that $\pi^{-1}(0)$ is a normal
crossings divisor with components $E_1,\dots,E_N$, let $Z$ be a
component of $E_I:=\bigcap_{i\in I} E_i$ for some family
$I\subset\{1,\dots,N\}$, and let $w_i>0$ be given real numbers for
$i\in I$. If $p$ is some point of $Z$, the local equations $z_i$,
$i\in I$ of the $E_i$ at $p$ can be completed to a local coordinate
system $z_1,\dots,z_n$ on $X_\pi$ at $p$. 
Then there exists a unique quasi-monomial
valuation on $X$ which is monomial on $X_\pi$ in the coordinates
$z_1,\dots,z_n$ with values $\nu(z_i)=w_i$ for $i\in I$ and
$\nu(z_i)=0$ otherwise. It is easily seen that $\nu$ does not depend
on the choice of $p\in Z$ nor on the way the $z_i$, $i\in I$ are
completed to a coordinate system at $p$. We will refer to $\nu$ as the
monomial valuation on $X_\pi$ with center $Z$ and value $\nu(E_i)=w_i$
on $E_i$.

Conversely, it follows from Hironaka's principalization theorem that
every quasi-monomial valuation centered at 0 is obtained in
that way.

\begin{Def}
  We denote by $\hcVdiv\subset\hcV$ (resp.\ $\hcVqm\subset\hcV$)
  the set of divisorial (resp.\ quasi-monomial)
  valuations centered at 0 in $X$, and by 
  $\cVdiv\subset\cV$, $\cVqm\subset\cV$ 
  the corresponding sets of normalized valuations.
\end{Def}
Note that every quasi-monomial valuation has home $X$ itself, and in particular the trivial valuation $\nu_0$ is definitely not a quasi-monomial valuation. Thus $\hcVqm$ (resp. $\hcVdiv$) is just the open cone over $\cVqm$ (resp. $\cVdiv$). 

The vanishing order at 0 $\nu_\fm(f) \= \max\{ k\in \N^*, f \in \fm^k\}$ is a
valuation that plays a special role in $\cV$. It is the divisorial valuation associated to the blowup of 0, and is also the monomial valuation on $X$ 
(in any local coordinates at 0)
with all weights equal to $1$. It is normalized, and satisfies the
following minimality property: for all $\nu \in \cV$ and all $f \in R$
we have $\nu(f)\ge \nu_\fm(f)$.
%
%
\subsection{Simple modifications and dual complexes}
We now aim at shedding some light on the geometry of $\cV$. 
Our presentation is inspired by~\cite{KoSo},~\cite{KKMS} 
and~\cite{thuillier2}.

A modification $\pi\in\fB$ is \emph{simple} if its
exceptional divisor $\pi^{-1}(0)=E_1+\dots+E_N$ has simple normal
crossings, and is such that $E_I:=\bigcap_{i\in I}E_i$ is irreducible
(or empty) for each $I\subset\{1,\dots,N\}$. Modifications with simple
normal crossing exceptional divisor are cofinite in $\fB$ as a
consequence of Hironaka's principalization theorem, and the
irreducibility condition can then be reached by further blowing-up
components of intersections of exceptional primes sufficiently many
times. Thus simple modifications are cofinite in $\fB$.

To a simple modification $\pi$ is attached a simplicial complex
which encodes the incidence relations between the $E_i$'s. 
Its vertices correspond to the $E_i$, and
its faces to subsets $I\subset\{1,\dots,N\}$ such that $E_I$ is
non-empty. We shall now show how to embed 
the underlying space of this complex into 
both $\R^N$ and into $\cV$.

First define $\hat{\Delta}(\pi)\subset\R^N$ to be the union of
all rational simplicial cones 
$\hat{\sigma}_I:=\sum_{i\in I}\R_+e_i\subset\R^I$, where $I$ runs over all subsets of
$\{1,\dots,N\}$ with $E_I\ne\emptyset$, and $(e_1,\dots,e_N)$ is the
canonical basis of $\R^N$.

Now pick a point $w\in\R^N$ lying in
$\hat{\Delta}(\pi)\setminus\{0\}$, 
and let $I=\{i, w_i>0\}$, so
that $\hat{\sigma}_I$ is the maximal cone to which $w$ belongs.  Since
$E_I$ is non-empty and irreducible, we can consider the quasi-monomial
valuation $\nu_w$ on $X$ that is monomial on $X_\pi$ with center $E_I$
and value $\nu_w(E_i)=w_i$ for each $i\in I$, 
as explained in the last section.  One easily verifies
\begin{Lemma}
  The map $w \mapsto \nu_w$ yields an embedding of
  $\hat{\Delta}(\pi)\setminus\{0\}$ into the set $\hcVqm$ 
  of quasi-monomial valuations in $\hcV$.
\end{Lemma}
Let us now describe the effect of our choice of normalization.  A
point $w\in\hat{\Delta}(\pi)$ induces a valuation $\nu_w$ belonging to
$\cV$ iff $\nu_w(\fm)=\sum_i b_i w_i = 1$, with
$b_i:=\ord_{E_i}(\fm)$, a positive integer.  With this in mind we
define $\Delta(\pi)$ as the trace of $\hat\Delta(\pi)$ on the integral
affine hyperplane $\sum_i b_iw_i=1$. It is a simplicial complex with
simplicial faces cut out by integral affine equations, all of the form
$\sigma_I=\hat{\sigma}_I\cap\Delta(\pi)$ for some $I$ with $E_I$
non-empty. The above lemma shows that $\Delta(\pi)$ embeds as a
topological subspace of $\cV$, which will be denoted by $|\Delta(\pi)|$.
\smallskip

We can endow the topological space $|\Delta(\pi)|$ with the 
structure of \emph{integral simplicial complex} on $|\Delta(\pi)|$.
We refer to~\cite[pp.69-70]{KKMS} for a general definition.
Suffice it to say that $|\Delta(\pi)|$ comes equipped
with faces $|\sigma_I|$, defined as the homeomorphic images 
of the faces $\sigma_I$ of $\Delta(\pi)$ 
by the above embedding into $\cVqm$, 
and by a free $\Z$-module 
$\Aff(\pi)_\Z\subset C^0(|\Delta(\pi)|)$
of rank $N$, induced by the 
restrictions to $\Delta(\pi)$ of all integral 
linear functions on $\R^N$.

The data $(|\Delta(\pi)|,\Aff(\pi)_\Z)$ enables us to 
completely recover the simplicial complex $\Delta(\pi)\subset\R^N$, 
up to unimodular transformations. Indeed, an appropriate choice of a $\Z$-basis $f_1,\dots,f_N$ for $\Aff(\pi)_\Z$ defines a map $|\Delta(\pi)|\to\R^N$ whose image is exactly $\Delta(\pi)$. We thus think of the data $(|\Delta(\pi)|,\Aff(\pi)_\Z)$ as $\Delta(\pi)$ itself, with its faces and its integral affine functions, but embedded in $\cVqm$ instead of an Euclidean space.
\begin{Def}
  The \emph{dual (simplicial) complex} $\Delta(\pi)$ of a simple modification
  $\pi$ is the integral simplicial complex consisting of the
  topological space $|\Delta(\pi)|\subset\cVqm$ endowed with the
  lattice $\Aff(\pi)_\Z$ of integral affine functions.
\end{Def}
We introduce the following useful terminology.
\begin{Def}\label{def:monoblw}
 If $\pi$ is a simple modification with exceptional divisors $E_i$,
and $\mu$ is the blow-up of $X_\pi$ with center $E_I$ for some $I$,
then $\mu$ is called a \emph{monomial blow-up}.  Two simple
modifications $\pi$ and $\pi'$ are said to be \emph{compatible} when
one is obtained from the other by a sequence of monomial blow-ups.
\end{Def}
The subset $|\Delta(\pi)|\subset\cVqm$ is not uniquely
determined by $\pi$. If $\pi$ is a simple modification, and $\mu$ is
the monomial blow-up of $X_\pi$ with center $E_I$ for some $I$, then
$\Delta(\pi\circ\mu)$ is obtained from $\Delta(\pi)$ by a barycentric
subdivision of the face $\sigma_I$. Hence
$|\Delta(\pi\circ\mu)|=|\Delta(\pi)|$ as subsets of 
$\cVqm$.\footnote{Note, however, that 
$\Aff(\pi)\subsetneq \Aff(\pi\circ \mu)$.} This is
essentially the only case: one can show that given two simple
modifications $\pi$ and $\pi'$, $|\Delta(\pi)|=|\Delta(\pi')|$ holds
iff there exists another simple modification $\pi''$ which is
compatible with both $\pi$ and $\pi'$. The dual complex $\Delta(\pi')$ of $\pi'$ compatible with $\pi$ is thus obtained by a sequence of barycentric subdivisions of $\Delta(\pi)$. A function on $|\Delta(\pi)|$
will be said to be \emph{piecewise integral affine} if it belongs to
$\Aff(\pi')_\Z$ for some $\pi'$ compatible with $\pi$
\begin{Def}
  A \emph{simplex} $|\sigma|$ of $\cVqm$ is a face of 
  $|\Delta(\pi)|$ for a simple modification $\pi$. 
  It is endowed with an integral affine structure
  $\Aff(|\sigma|)_\Z\subset C^0(|\sigma|)$, 
  defined as the restriction of $\Aff(\pi)_\Z$ to $|\sigma|$.
\end{Def}
The integral affine structure on a simplex $|\sigma|$ is independent of
the choice of $\pi$ above, and functions in $\Aff(|\sigma|)_\Z$
will be said to be integral affine on $|\sigma|$. 
If $\pi$ is a simple modification such that $|\sigma|$ is
a face of $|\Delta(\pi)|$, then, with notation as above,
the elements of $|\sigma|$ can be described as the
monomial valuations $\nu_w$ on $X_\pi$ centered on $E_I$ taking
values $w_i$ on $E_i$ for $i\in I$ and 
such that $\sum_{i\in I}b_iw_i=1$. 
The coordinates $t_i:=b_i w_i$ identify $|\sigma|$ with the
standard simplex in $\R^I$, and $\Aff(|\sigma|)_\Z$ corresponds to the
standard integral functions $t\mapsto\sum_i a_i t_i$, $a_i\in\Z$.  
The set $\cVdiv\cap|\sigma|$ 
of divisorial valuations contained in $|\sigma|$
is exactly the set of points with rational $w$-coordinates. 
In particular, it is dense in $|\sigma|$.

If $|\sigma|$, $|\sigma'|$ are two simplices in the same $|\Delta(\pi)|$, 
then $|\sigma|\cap|\sigma'|$ is
either empty or a simplex having integral affine structure
$\Aff(|\sigma|\cap|\sigma'|)_\Z$ equal to the restriction of both
$\Aff(|\sigma|)_\Z$ and $\Aff(|\sigma'|)_\Z$ to $|\sigma|\cap|\sigma'|$.

The main relation between the affine structure on simplices and formal
functions in $\hat R$ is the following.
\begin{Lemma}\label{lem:meroaffine}
  If $f\in\hat R$ is a formal function, the restriction of $\log|f|$
  to any simplex $|\sigma|$ is convex, continuous,
  and locally piecewise affine.
\end{Lemma}
\begin{proof}
  Let $\pi$ be a simple modification such that $|\sigma|=|\sigma_I|$ 
  occurs as a face of $|\Delta(\pi)|$ as above. If $p$ is a point of
  $E_I$, and local equations $z_i$ of the $E_i$ at $p$ have been
  completed to local coordinates $z_1,\dots,z_n$ on $X_\pi$ at $p$, 
  we can write 
  $f\circ\pi=\sum_\alpha a_\alpha z^{\alpha}$. For any $w\in\sigma_I$,
  with corresponding valuation $\nu_w\in|\sigma|$, we then have
  $\nu_w(f)=\min\{w\cdot\alpha, a_\alpha\neq 0\}$, so that
  $\log|f|(\nu_w)= -\nu_w(f)$ is the maximum of the family of
  \emph{integral} affine functions $w\mapsto -w\cdot\alpha$ such that
  $a_\alpha\neq 0$. The maximum can be taken over finitely many
  of these functions and then all the conclusions are immediate.
\end{proof}
%
%
\subsection{Retractions and piecewise affine functions}\label{sec:affine}
Let $\Delta(\pi)$ be the dual complex of a simple
modification $\pi$, with exceptional primes $E_1,\dots,E_N$. Given any
valuation $\nu \in \cV$, let $I$ be the family of primes $E_i$ that
contain the center $Z$ of $\nu$ on $X_\pi$, \ie such that
$\nu(E_i)>0$. There exists a unique valuation $r_\pi(\nu)$ which is
monomial on $X_\pi$, with center $E_I$ and such that
$r_\pi(\nu)(E_i)=\nu(E_i)$ for all $i\in I$.

We thus get a continuous retraction map $r_\pi: \cV \to|\Delta(\pi)|$,
which in fact only depends on the subset $|\Delta(\pi)|\subset\cV$. 
Indeed, it is
enough to check that $r_\pi=r_{\pi'}$ whenever $\pi'=\pi\circ\mu$ for a
monomial blow-up $\mu$ of $X_\pi$, and this is easily done.
\begin{Prop}
  For any simple modification $\pi$, there exists a
  naturally defined retraction map
  $r_\pi:\cV\to|\Delta(\pi)|$. Furthermore, if $\pi'$ is any other
  simple modification dominating $\pi$, the restriction
  $r_{\pi',\pi}:|\Delta(\pi')|\to|\Delta(\pi)|$ of $r_\pi$ to
  $|\Delta(\pi')|$ is integral affine, in the sense that
  $r_{\pi',\pi}^{*}$ sends $\Aff(\pi)_\Z$ into $\Aff(\pi')_\Z$.
\end{Prop}
Using $r_\pi$, we can identify $\Aff(\pi)_\Z$ with the subgroup
$r_\pi^*\Aff(\pi)_\Z\subset C^0(\cV)$.
\begin{Def} 
  A \emph{piecewise affine} function 
  on $\cV$ is an element of the union $\PA(\cV)$ of all
  $\Aff(\pi)_\Z\otimes \R\subset C^0(\cV)$. 
  It is said to be \emph{integral}
  when it belongs to some $\Aff(\pi)_\Z$.
\end{Def}

The inductive system of all simple modifications $\pi$ gives rise to
an inductive system of integral simplicial complexes $|\Delta(\pi)|$ with integral affine maps
$r_{\pi',\pi}:|\Delta(\pi')|\to|\Delta(\pi)|$.
The projective limit $\varprojlim_\pi|\Delta(\pi)|$ of compact
topological spaces is compact, and there is a natural map $\cV \to
\varprojlim_\pi|\Delta(\pi)|$ given by the collection $r_\pi$. This
map is injective because, for any two distinct valuations $\nu_1,
\nu_2\in\cV$ there exists a simple modification $\pi\in \fB$ such that
the centers of $\nu_1$ and $\nu_2$ on $X_\pi$ are distinct.

On the other hand, we have seen that any quasi-monomial valuation lies
in some $|\Delta(\pi)|$. Whenever $\pi'$ dominates $\pi$ the inclusion
$|\Delta(\pi)|\subset|\Delta(\pi')|$ realizes $\Delta(\pi)$ as a
integral simplicial sub-complex of $\Delta(\pi')$. We thus have injections
\[
\varinjlim_\pi|\Delta(\pi)|=\cVqm\hookrightarrow
\cV\hookrightarrow\varprojlim_\pi|\Delta(\pi)|.
\]
The first space being dense in the last, and
$\cV$ being compact we conclude that 
$\cV=\varprojlim|\Delta(\pi)|$.  We have thus proved
\begin{Thm}\label{thm:strucV}
  The space $\cV$ with its piecewise integral affine functions
  $\PA(\cV)_\Z$ can be naturally identified with the projective limit
  of the system of integral simplicial complexes
  $r_{\pi',\pi}:|\Delta(\pi')|\to|\Delta(\pi)|$.
  
  The space $\cVqm$ of quasi-monomial valuations can be naturally
  identified with the inductive limit $\varinjlim_\pi |\Delta(\pi)|$.
 \end{Thm}
In view of the preceding result, it is natural to introduce the
\emph{strong topology} on $\cVqm$ as the inductive limit
topology.  Since each dual complex is a finite union of simplices,
a function $g$ on $\cVqm$ is continuous iff its restriction to every
simplex $|\sigma|$ is continuous.
%
%
\subsection{Ideals, Cartier divisors and piecewise affine functions}
We now relate piecewise affine functions on $\cV$ to primary ideals
of $R$ and Cartier divisors on $\fX$.

First note that a Weil divisor $Z\in \div(\fX)$ induces a real-valued
function $g_Z:\cVdiv \to \R$ by setting
$g_Z(\nu)=\ord_E(Z)/\ord_E(\fm)$ whenever $\nu$ is the normalized
valuation proportional to a divisorial valuation $\ord_E$. This
induces a natural bijection between Weil divisors on $\fX$ 
and functions on $\cVdiv$. 
\begin{Prop}\label{p:cartier}
  If $Z$ is a Cartier divisor on $\fX$, the function $g_Z:\cVdiv\to\R$
  extends by continuity to a function in $\PA(\cV)$, and this induces
  an isomorphism $\Cdiv(\fX)\to \PA(\cV)$. If $\fa$ is a primary ideal
  of $R$, the Cartier divisor $Z(\fa)$ it determines is associated to
  the function $\log|\fa|$, and $\PA(\cV)_\Z$ coincides with the
  subgroup of $C^0(\cV)$ generated by $\log|\fa|$ for $\fa$ ranging
  over all primary ideals of $R$.
\end{Prop}
\begin{proof}
  If $\pi$ is a simple modification, the vertices of $\Delta(\pi)$
  correspond to the exceptional primes of $\pi$. This yields an
  identification $\div(\pi)_\Z\to\PA(\pi)_\Z$, by sending a divisor
  $Z_\pi$ to the unique function on $\Delta(\pi)$ that is affine on
  each face and whose value at each vertex $E$ is given by
  $\ord_E(Z_\pi)/\ord_E(\fm)$. Now, if $\pi'=\pi\circ\mu$ is another
  simple modification, it is easily checked that
  $\mu^{*}:\div(\pi)\to\div(\pi')$ corresponds to the map
  $(r_{\pi,\pi'})^{*}:\PA(\pi)\to \PA(\pi')$. This proves that the
  above isomorphisms induce an isomorphism between
  $\Cdiv(\fX)_\Z=\varinjlim\div(\pi)_\Z$ and $
  \PA(\cV)_\Z=\varinjlim\PA(\pi)_\Z$.
  
  Now a given $\Z$-divisor $D$ on some $X_\pi$ can be written as a
  difference of two $\Z$-divisors $D_1-D_2$, where each $D_i$ is
  generated by its global sections on $X_\pi$. This means that
  $\cO(D_i)=\pi^{-1}\fa_i$ for some primary ideals $\fa_i$ in $R$, and
  this shows that the group generated by functions of the form
  $\log|\fa|$ with $\fa$ a primary ideal coincides with $\PA(\cV)_\Z$.
\end{proof}
\begin{Prop}\label{p:dense} 
  The subspace $\PA(\cV)$ is 
  dense in $C^0(\cV)$, endowed with the topology of uniform convergence.
\end{Prop}
\begin{proof}
  This is a straightforward consequence of the so-called lattice version 
  of the   Stone-Weierstrass
  theorem. 
  Indeed, we have $|\fa+\fb|=\max(|\fa|,|\fb|)$ as functions on $\cV$ 
  for any   two primary ideals $\fa$, $\fb$, and it follows immediately that   the closure of the
 $\Q$-vector subspace of $C^0(\cV)$ spanned by all functions of the form $\log|\fa|$ is a closed 
$\R$-subspace of $C^0(\cV)$ stable by taking max. Since it also separates points of $\cV$ in a tautological way, the lattice version of the Stone-Weierstrass theorem implies that this space is equal to $C^0(\cV)$. But it also coincides with the closure of $\PA(\cV)$, and the result follows.
\end{proof}
%
%
\subsection{Contractibility}
We conclude with the following result due to Berkovich~\cite[p.116]{ber}. It
is also a special case of a more general result of \cite{thuillier2}.
\begin{Thm}[\cite{ber,thuillier2}] 
The spaces $\cV$, $\cVqm$ are contractible.
\end{Thm} 
\begin{proof} Let $x_1,\dots,x_n$ be coordinates on $X$ at $0$,
and denote by $\cVm\subset \cV$ the space of all monomial valuations,
\ie determined by their values on $x_i$. Clearly $\cVm$ can be
identified with $\{(s_i)\in[1,+\infty]^n,\, \min s_i = +1 \}$, and is
thus contractible.  Now define the retraction $r: \cV\to \cVm$ by:
\[r(\nu)(f)=\min_\alpha \{\nu(a_\alpha x^{\alpha})\}=
\min\{\nu(x^{\alpha}), a_\alpha \neq0 \},\]
where $f=\sum_\alpha a_\alpha x^{\alpha}$ in multi-index notation. 
We construct a homotopy
$h_t$ from $r$ to the identity as follows. For each multi-index
$\beta\in\N^n$, let $T_\beta$ be the truncation operator on
$k[x_1,\dots,x_n]$ given by $T_\beta(\sum_\alpha a_\alpha
x^{\alpha})=\sum_{\alpha\geq\beta}a_\alpha x^{\alpha}$.  Here $\alpha
= (\alpha_i) \ge (\beta_i)$ iff $\alpha_i \ge \beta_i$ for each
$i$. Then for $t\in[0,1]$, we set
\[h_t(\nu)(f)\=\min_\beta \{\nu(T_\beta(f))- |\beta| \log t\}.\]
Note that the minimum is attained for only finitely many $\beta$'s.
To see that the function $h_t(\nu)$ is a valuation, it is enough to check 
its additivity on the product of two functions.
Pick $f,g\in R$, and let $\alpha$ (resp. $\beta$) be the 
minimal index for the lexicographic order on $\N^n$ 
such that $h_t(\nu)(f) = \nu(T_\alpha(f))- |\alpha| \log t$
(resp. $h_t(\nu)(g) = \nu(T_\beta(g))- |\beta| \log t$).
Set $\gamma = \alpha + \beta$. Then  
$\nu(T_{\gamma}(fg)) = \nu( \sum_{\alpha'+\beta' = \gamma} T_{\alpha'}(f) T_{\beta'}(g)) = 
\nu(T_\alpha(f)) + \nu(T_\beta(g))$,  hence $h_t(\nu) (fg) \le 
h_t(\nu) (f)+ h_t(\nu) (g)$. The reverse inequality is clear.

It
is easy to see that $h_0$ is the identity on $\cV$, and $h_1$ coincides with
$r$ in view of the relation $\min_\beta \nu(\sum_{\alpha\geq\beta}
a_\alpha x^{\alpha})=\min_\beta \nu(a_\beta x^{\beta})$. Continuity of
$h_t(\nu)$ in $(t,\nu)$ is also easily checked, and $h_t$ is the
identity on $\cVm$ for each $t$.
\end{proof}
For each simple modification $\pi$, the retraction
$r_\pi:\cV\to|\Delta(\pi)|$ 
can also be deformed to the identity map of
$\cV$ relatively to $|\Delta(\pi)|$. 
This is done in \cite{thuillier2}
using the language of Berkovich spaces, relying on the fact that the
normal crossing exceptional divisor of $\pi$ endows $X_\pi$ with a
toroidal structure.
\begin{Remark}
The valuation spaces $\cVqm$ and $\cV$ in dimension $n=2$ were studied
by the last two authors in~\cite{valtree}. In this case, $\cV$ is
naturally a tree---\emph{the valuative tree}--- and the piecewise
affine integral structure corresponds to the parameterization of $\cV$
by skewness.  Moreover, $\cVqm$ is a dense subtree of $\cV$ obtained
by removing all the ends of $\cV$.  This makes it possible to
understand the valuations in $\cV\setminus\cVqm$.  In higher
dimensions, the structure of $\cV\setminus\cVqm$ is much more
complicated. For instance, it contains the valuation space $\cV_Y$ for
each smooth (formal) subvariety $Y\subset X$ through $0$. In dimension
two, $Y$ is a germ of a curve, and $\cV_Y$ is then a singleton, called a
curve valuation.
\end{Remark}

%
%
%
%
\section{Positivity and convexity}\label{sec:pos}
%
%
\subsection{Nef Weil divisors}\label{sec:pos1}
Recall that to any ideal $\fa\subset \hat{R}$ is associated a Weil
divisor $Z(\fa)\le 0$ which is Cartier iff $\fa$ is primary.  For any
two ideals $\fa$ and $\fb$ in $\hat R$, we have
\begin{align}
  Z(\fa \cdot \fb) &= Z(\fa)+ Z(\fb);\label{e1} \\
  Z(\fa + \fb) &= \max\{ Z(\fa), Z(\fb)\}. \label{e2}
\end{align}
Here the max is taken component-wise. This implies that the set of all
$Z(\fa)$'s is a semi-group preserved by taking $\max$. It is the
smallest such semi-group containing the divisors $Z(f)$ determined by
all $f\in \hat R$, since $Z(\fa)=\max_k Z(f_k)$ whenever the $f_k$'s
generate $\fa$.
\begin{Def}
  The cone of \emph{nef Weil divisors} 
  in $\div(\fX)$ is defined as the 
  closed convex cone in $\div(\fX)$ generated by the  divisors $Z(\fa)$,
  where $\fa$ runs over all ideals of $\hat R$.
\end{Def}
Note that a nef Weil divisor $Z$ satisfies $Z\leq 0$, by continuity.
If $\fa$ is an arbitrary ideal, then we have
$Z(\fa)=\lim_{k\to\infty}Z(\fa+\fm^k)$, so that the $Z(\fa)$'s with
$\fa$ a primary ideal of $R$ already generate the nef cone. 
Since the $Z(\fa)$'s make up a semi-group preserved by taking $\max$,
the nef cone of $\div(\fX)$ is also preserved by the same operation.

In the special case of Cartier divisors, we have the following
characterization, which explains the terminology.
\begin{Prop}\label{prop:nefcartier}
  A Cartier divisor $Z$ on $\fX$ is nef (as a Weil divisor)
  iff $Z_\pi$ is 
  $\pi$-nef on $X_\pi$ for some (hence any) determination $\pi$ of $Z$,
  that is, $Z_\pi$ has non-negative degree on any compact curve 
  contained in $\pi^{-1}(0)$.
\end{Prop}
\begin{Lemma}
  Let $Z$ be a nef Weil divisor on $\fX$, and let 
  $Z_\pi$ be its incarnation on some
  $X_\pi$. For any exceptional prime $E$ of $\pi$, the restriction
  $Z_\pi|E$ of the numerical class of $Z_\pi$ to $E$ is pseudo-effective,
  \ie it is a limit of classes of effective $\Q$-divisors on $E$. 
\end{Lemma}
\begin{proof}
  Since the set of pseudo-effective classes of $E$ 
  is a closed convex cone, it is enough to check this 
  when $Z=Z(\fa)$ for some primary ideal $\fa$. 
 But $Z_\pi$ is the divisorial part of $\pi^{-1}\fa$, 
  thus the base-locus of the global sections of $\cO(Z_\pi)$ 
  on $X_\pi$ has codimension at least 2. 
  It follows that $E$ cannot be contained in this base-locus, 
  so $\cO_E(Z_\pi)$ is a line bundle with at least one non-zero 
  global section. In particular, it is pseudo-effective. 
\end{proof}
\begin{proof}[Proof of Proposition~\ref{prop:nefcartier}]
  Let $\pi$ be a determination of $Z$. If $Z_\pi$ is not $\pi$-nef 
  on $X_\pi$, 
  there  exists a curve $C$ in $\pi^{-1}(0)$ such that the numerical 
  class $Z_\pi|C$ is not pseudo-effective. 
  Now pick $\pi'=\pi\circ\mu\in\fB$ dominating the blow-up of $C$ 
  in $X_\pi$, and let $E$ be an exceptional prime of 
  $\pi'$ dominating $C$. Then the restriction of $\mu^{*}Z_\pi=Z_{\pi'}$ 
  to $E$, which is the same as $(\mu|E)^{*}(Z_\pi|C)$, 
  cannot be pseudo-effective since $Z_\pi|C$ is not.
  This leads to a contradiction in view of the lemma, since $Z_{\pi'}$ 
  is the incarnation of the nef divisor $Z$ on $X_{\pi'}$. 
\end{proof} 
We record the following monotonicity property of nef divisors.
\begin{Prop}\label{p:decrease}
  Let $Z$ be a nef Weil divisor on $\fX$, 
  and let $Z_\pi\in\div(\fX)$ be the incarnation of $Z$ on $X_\pi$,
  viewed as a Cartier divisor on $\fX$.
  Then $Z\le Z_\pi$. 

  In particular, any nef Weil divisor $Z\ne0$ satisfies
  $Z\le c Z(\fm)$ for some $c>0$.
\end{Prop}
\begin{proof}
  Since $Z\mapsto Z_\pi$ is linear and continuous
  we can assume $Z=Z(\fa)$ for some primary ideal $\fa$.
  Let $\pi'=\pi\circ\mu$ be a determination of $Z$ and
  write $L:=Z_{\pi'}$. 
  Then $\cO(L)$ is generated by its global sections
  with respect to $\pi'$, thus \emph{a fortiori} with respect to
  $\mu$. In other words, we have
  $\cO(L)=\mu^{*}\mu_{*}\cO(L)$. Now the sheaf of ideals
  $\cO(\mu_{*}L)$ contains $\mu_{*}\cO(L)$ since $L\leq 0$, 
  so that $\cO(L)\subset\cO(\mu^{*}\mu_{*}L)$.
  Since $L$ is exceptional, 
  this implies $Z=Z_{\pi'}\le Z_\pi$.
\end{proof}
The proof also shows that if $Z$ is any nef Weil divisor
on $\fX$, then $Z_{\pi'}\le Z_\pi$ whenever $\pi'\ge\pi$, 
and $Z=\inf_\pi Z_\pi$.
%
%
\subsection{Formal psh functions}\label{sec:convx}
We now interpret the previous positivity notion 
at the level of functions on the valuation space $\cV$. 
The main technical tool required to study 
them---multiplier ideals---will be introduced in the next section.
\begin{Def}\label{def:convex}
  A function $g : \cVqm \to \R$ is called a \emph{formal psh function} 
  iff it is continuous
  in the strong topology and if the Weil divisor $Z(g)$ 
  on $\fX$ determined by $g|\cVdiv$ is nef.
\end{Def}
As should be expected from the name, formal psh functions 
share many properties with usual the psh functions on $\C^n$.
In particular, $\log|f|$ is a formal psh function for any 
formal function $f\in\hat{R}$, see Lemma~\ref{lem:meroaffine}.
\begin{Remark}
We will prove later on that any nef Weil divisor $W$ on $\fX$
can conversely be written (in a unique way) as $W=Z(g)$ for some formal psh function $g$ (see Proposition~\ref{p:extens}), so that nef Weil divisors on $\fX$ and formal psh functions on $\cV$ are
two manifestations of the same object. Let us already here point out
what the difficulty is. Given a simplex in $|\sigma|$ in $\cV$
we must show that the function $g:|\sigma|\cap\cVdiv\to\R$
defined by $W$ extends to a continuous, convex function on $|\sigma|$.
Suppose $|\sigma|$ is one-dimensional, for simplicity, so
that $|\sigma|$ can be identified with the interval $[0,1]$ and
$|\sigma|\cap\cVdiv$ with $[0,1|\cap\Q$.
The nefness assumption on $W$ implies that 
$g:[0,1]\cap\Q\to\R$ is convex.
However, while g then
   automatically extends as a continuous function on the open
   interval $(0,1)$, there may not be a continuous extension to
   the \emph{closed} interval $[0,1]$.
An example is
given by $g(x)=-x$ for $0\le x<1$ and $g(1)=0$.
\end{Remark}  
If $Z$ is a Weil divisor on $\fX$ and $\pi$ is a simple modification,
recall that we denoted by $Z_\pi$ the Cartier divisor in $\fX$ 
determined by  the incarnation of $Z$ in $X_\pi$.
In terms of the associated function $g_Z$ on
$\cVdiv$, we see that $Z_\pi$ corresponds to the unique affine
function in $\Aff(\pi)$ that coincides with $g_Z$ on the vertices of
$\Delta(\pi)$, \ie to the piecewise affine interpolation of $g_Z$ on
$\Delta(\pi)$.  Now suppose that $g$ is a formal psh function. Then by
Proposition~\ref{p:decrease}, $Z(g)\le Z(g)_\pi$ for any $\pi$. Both
functions $g_{Z_\pi}$ and $g\circ r_\pi$ belong to 
$r_\pi^{*}C^0(|\Delta(\pi)|)$, but they do not coincide, since the
restriction of $g$ to $|\Delta(\pi)|$ is not piecewise affine in
general. Using the notion of compatible blowups from
Definition~\ref{def:monoblw}, we have:
\[
g\circ r_\pi= \inf\, \{ g_{Z_{\pi'}},\, \pi'\
\text{compatible with $\pi$}\}.
\]
Using this fact, we conclude
\begin{Prop}\label{prop:fpshdecreasing}
  For any dual complex $\Delta(\pi)$ and any
  formal psh function $g$, we have $g \le g \circ r_\pi$.
\end{Prop}
In dimension 2, $Z_\pi$ is nef whenever $Z$ is, and
thus $g\circ r_\pi$ is formal psh when $g$ is. 
These two properties fail in the higher dimensional case.

If $\pi'\geq\pi$ are two simple modifications, we have
$|\Delta(\pi)|\subset|\Delta(\pi')|$, thus $r_\pi\circ
r_{\pi'}=r_\pi$. Hence $g\circ r_\pi$ decreases as $\pi$
increases. Using this, we extend $g$ to $\cV\setminus\cVqm$ as follows:
\begin{Def}
  Given a formal psh function $g:\cVqm\to\R$, we extend 
  $g$ to $\cV$ by setting
  \[
  g(\nu) \= \inf_\pi g \circ r_\pi(\nu)
  \]
  for any valuation $\nu\in\cV$.
\end{Def}
Since the restriction of $g$ to any dual complex $|\Delta(\pi)|$ is
continuous, we see that $g\circ r_\pi$ is continuous in the weak
topology of $\cV$ for each $\pi$, and thus the extension
$g:\cV\to[-\infty,0]$ is \emph{upper} semicontinuous, 
as in the case of the usual psh functions.
%
%
%
%
\section{Valuative multiplier ideals}\label{sec:mult}
The goal of this section is to attach to each nef Weil divisor $W$ on
$\fX$ a valuative multiplier ideal $\elde(W)\subset\hat R$, in such a
way that $\elde(cZ(\fa))$ coincides with the usual multiplier ideal
$\cJ(\fa^c)$ when $\fa$ is a primary ideal of $R$ and $c>0$. We will
then show that the nef Weil divisors $\frac1k Z(\elde(kW))$ approximate in a very precise way 
a given nef Weil divisor $W$ from above.
This will enable us to further investigate formal psh functions.  
In particular, it will turn out that formal psh functions and
nef Weil divisors are equivalent notions.
%
%
\subsection{Thinness}
Consider a modification $\pi:X_\pi\to X$ in $\fB$.
Denote by $K_{\pi}$ the relative canonical divisor of $\pi$. 
It is the effective divisor on $X_\pi$ determined by the
Jacobian determinant $J{\pi}$ of $\pi$. 
When $\pi'\ge\pi$, $\pi' = \pi\circ \mu$, the chain rule 
yields $K_{\pi'}=\mu^*K_{\pi}+K_{\mu}$, where $K_{\mu}$, 
the relative canonical divisor of $\mu$, is $\mu$-exceptional. 
Thus there exists a unique Weil divisor $K$ on $\fX$, 
the \emph{canonical divisor},
whose incarnation on each $X_\pi$ is $K_\pi$.
For our purposes, a different but closely related Weil 
divisor will be important.
\begin{Def}
  The \emph{thin divisor} $A$ is the
  Weil divisor on $\fX$ whose incarnation in $X_\pi$ 
  equals $A_\pi\=K_\pi+\sum E$ 
  where $E$ ranges over all exceptional primes of $\pi$.
  The corresponding function $A:\cVdiv\to\R_+$ is called
  \emph{thinness}.
\end{Def}
Let $|\sigma|$ be a simplex in $\cVqm$, determined by a
collection $E_i$, $i\in I$ of 
exceptional primes of a simple modification $\pi$.
If $\nu\in|\sigma|$ is divisorial, a standard computation 
(see for instance~\cite{howald}) shows that 
$A(\nu)=\sum_{i\in I}w_iA(\ord_{E_i})$ with $w_i=\nu(E_i)$. 
This yields:
\begin{Prop}
  The thinness function $A:\cVdiv\to\R_+$ extends uniquely
  to a strongly continuous function $A:\cVqm \to \R_+$,
  which is integral affine on each simplex.
\end{Prop}
The thinness function satisfies the following monotonicity property.
\begin{Prop}
  For any simple modification $\pi$, we have $A\geq A\circ r_\pi$ on
  $\cVqm$, and the inequality is strict 
  outside the dual complex $|\Delta(\pi)|$.
\end{Prop}
The proof is a standard computation, \cf~\cite[Lemma~3.11]{kollar}.

In particular we have $A\geq A(\nu_\fm)=\dim X$ on $\cVqm$.
We also see that $A\circ r_\pi:\cVqm\to\R$ is increasing 
with respect to $\pi$, so we can extend $A$ to all of $\cV$ as follows. 
\begin{Def} 
  The thinness function $A:\cV\to[0,+\infty]$ is the lower
  semicontinuous function defined by $A=\sup_\pi A\circ r_\pi$.
\end{Def}
Beware that the Weil divisor $-A$ on $\fX$ is not nef, \ie $-A$ is 
not a formal psh function on $\cV$.  
%
%
\subsection{Multiplier ideals of ideals}\label{sec:def-mult}
We start by recalling the usual geometric 
definition of the multiplier ideals following~\cite[\S9.2]{laz2}.

Let $c>0$ and let $\fa\subset R$ be an ideal 
(not necessarily primary).
Pick a log-resolution $\pi:X'\to X$ of $\fa$.  Then $f\in R$
lies in the multiplier ideal $\cJ(\fa^c)$ iff
$\ord_E(f\circ\pi)+\ord_E(K_\pi)\ge[c\ord_E(\pi^{-1}\fa)]$
for\footnote{Here $[x]$ is the integral part (rounded down)
  of the real-number $x$.}
\emph{any} prime divisor $E$ of $X'$. If $\fa$ is primary, then $\pi$
can be taken to lie in $\fB$, and the above condition only needs to be
tested for the $\pi$-exceptional primes $E$. Using thinness, it can
then be rewritten as $\ord_E(f)+A(\ord_E)>c\ord_E(\fa)$.

The multiplier ideal $\cJ(\fa^c)$ satisfies the following fundamental
\emph{subadditivity} property, first established in~\cite{DEL}. 
We refer to~\cite[\S9.5]{laz2} for a proof.
\begin{Thm}\label{t:standard} 
  Let $\fa$ and $\fb$ be ideals in $R$.
  Then $\cJ(\fa^c\cdot\fb^{d})\subset\cJ(\fa^c)\cdot\cJ(\fb^{d})$
  for any $c,d>0$.
\end{Thm}
%
%
\subsection{Multiplier ideals of nef Weil divisors}
\begin{Def}
  Let $W$ be a nef Weil divisor on $\fX$.
  The \emph{multiplier ideal} $\elde(W)$ is defined as the
  stationary limit as $\e>0$ decreases to $0$ of the ideals
  \begin{equation}\label{e:valmul}
    \{f\in\hat R,\ Z(f) \le (1+\e)W +A\}.
  \end{equation}
  If $g$ is a formal psh function on $\cV$, then we define
  $\elde(g)\=\elde(Z(g))$.
\end{Def}
Note that we consider arbitrary formal functions $f\in\hat R$, not
just in $R$. Since $Z\leq 0$, the ideals in~\eqref{e:valmul} form an
increasing sequence as $\e >0$ decreases, hence are indeed stationary
by the Noetherian property. Since each of these ideals is integrally
closed, so is the multiplier ideal $\elde(W)$. Note however that it is
not primary in general. The definition is arranged so that
$\elde(W)=\elde((1+\e)W)$ for every $\e>0$ small enough. Also, it
coincides with the usual multiplier ideal when $W=cZ(\fa)$ for
some ideal $\fa$ of $R$:
\begin{Prop}\label{p:usual} 
  When $\fa\subset R$ is an ideal and $c>0$, 
  the ideal $\elde(c Z(\fa))\cap R$ coincides with the usual
  multiplier ideal $\cJ(\fa^c)$.
\end{Prop}
\begin{Remark} 
  Over the complex numbers, 
  if $f_1,\dots,f_k$ are generators of $\fa$, then $f$ is in the
  multiplier ideal of $\fa^c$ iff $|f|/\max_j|f_j|^c$ is locally
  square-integrable at 0.  
  This explains our notation $\elde$ above.
\end{Remark}
\begin{proof}[Proof of Proposition~\ref{p:usual}]
  This is essentially trivial when $\fa$ is primary, but the
  non-primary case requires some care since
  $\elde(c\log|\fa|)$ is defined by
  imposing conditions at valuations centered at $0$ on $X$ only. 
  
  Set $Z:=cZ(\fa)$ and $W=Z(\cJ(\fa^c))$. We know that
  $\ord_E(W)<\ord_E(Z)+A(\ord_E)$ for any 
  divisorial valuation $\ord_E$
  on $X$ (not necessarily centered at 0). 
  Indeed, $E$ can be realized as an exceptional
  prime of a log-resolution of $\fa$.
  As a consequence,
  $W<Z+A$ as Weil divisors on $\fX$. 
  Now, it is well known that $\cJ(\fa^c)=\cJ(\fa^{(1+\e)c})$
  for $0<\e\ll1$, so we get $W<(1+\e)Z+A$. 
  Thus $\cJ(\fa^c)\subset\elde(Z)\cap R$.

  To get the reverse inclusion, let $f\in\elde(Z)\cap R$, and pick a
  log-resolution $\pi:X'\to X$ of $\fa$. Let $\sum_i r_iE_i$ be the
  effective divisor defined by $\pi^{-1}\fa$. By definition, 
  $f$ belongs
  to $\cJ(\fa^c)$ iff $\ord_{E_i}(f)>cr_i-a_i$ for each prime $E_i$ of
  $\pi$, with $a_i:=1+\ord_{E_i}(J\pi)$.
  (Note that $\pi$ may not contract $E_i$ onto 0.)
  
  Let $E=E_i$ be one of these primes, pick a point $p$ of
  $E\cap\pi^{-1}(0)$ and let $z_1$ be a local equation of $E$ at $p$,
  completed to a local system of coordinates 
  $z_1,\dots,z_n$ on $X'$ at $p$. 
  If we let $\nu_w$ be the monomial valuation on $X'$ in the
  coordinates $z_1,\dots,z_n$ with weights $\nu(z_i)=w_i$, then
  $\ord_E=\nu_w$ for $w=(1,0,\dots,0)$, and the center of 
  $\nu_w$ on $X$
  is $\pi(p)=0$ if $w_i>0$ for all $i$. Since $f\in\elde(Z)$, 
  we deduce
  that $\nu_w(f)\geq(1+\e)c\nu_w(\fa)-A(\nu_w)$ for each $w$ with all
  $w_i>0$. But clearly $\nu_w(\fa)\to\ord_E(\fa)$ and
  $\nu_w(f)\to\ord_E(f)$ as $w>0$ converges to $(1,0,\dots,0)$, and
  $A(\nu_w)=w_1+\dots+w_n+\nu_w(J\pi)$ converges to $1+\ord_E(J\pi)$.
  We deduce that $\ord_E(f)\geq(1+\e)c\ord_E(\fa)-(1+\ord_E(J\pi))$
  for each $E$, and in particular $f\in\cJ(\fa^c)$ 
  as was to be shown.
\end{proof}

We will prove in Section~\ref{sec:pf} that the valuative multiplier
ideals satisfy the following fundamental properties, extending those
of the usual multiplier ideals.
\begin{Thm}\label{thm:approx}
  If $W$ is a nef Weil divisor on $\fX$, then
  \[
  W\leq Z(\elde(W)))\leq W+A.
  \]
\end{Thm}
\begin{Thm}\label{thm:subadditivity}
  If $W_1$ and $W_2$ are two nef Weil divisors on $\fX$, then
  \[
  \elde(W_1+W_2)\subset\elde(W_1)\cdot\elde(W_2).
  \]
  in $\hat R$.
\end{Thm}
%
%
\subsection{Consequences for formal psh functions}\label{sec:approx}
Using Theorems~\ref{thm:approx} and~\ref{thm:subadditivity} we
shall establish two important facts. 
The first is that nef Weil divisors and formal psh functions are 
equivalent notions, as already mentioned above.
\begin{Prop}\label{p:extens}
  For any nef Weil divisor $W$, there exists 
  a unique formal psh function
  $g:\cV\to[-\infty,0]$ such that $W=Z(g)$.
  Moreover, $g$ is continuous and convex on any simplex 
  in $\cV$.
\end{Prop}
\begin{proof}
  Let $g:\cVdiv\to\R$ be the function associated to $W$. 
  We have to show that $g$ can be continuously extended to $\cVqm$ 
  endowed with its strong topology. Concretely, this means that
  $g$ extends to a continuous function on each simplex $|\sigma|$
  in $\cV$. 

  The apparent difficulty is that $|\sigma|$ is a closed simplex.
  Indeed, the nefness of $W$ can be shown to imply that $g$ extends
  uniquely to a convex (and hence continuous!) function 
  on the \emph{interior} of $|\sigma|$. However, the behavior of
  $g$ at the boundary of $|\sigma|$ is harder to control, and
  that is why we use multiplier ideals.

  For each $k>0$, set $g_k:=\frac 1k\log|\elde(kW)|$. 
  This is a formal psh function, and $|g-g_k|\le A/k$ 
  on $\cVdiv$ by Theorem~\ref{thm:approx}. Since $A$ is 
  bounded on $|\sigma|$, $(g_k)_1^\infty$ is a 
  Cauchy sequence on $|\sigma|$,
  hence converges uniformly to a continuous function. 
  It is clear that the function thus obtained is the extension 
  of $g$ that we were looking for. The convexity of
  $g$ on simplices in $\cV$ is a consequence of 
  Proposition~\ref{prop:fpshdecreasing}.
\end{proof}
\begin{Cor}
  If $g:\cVqm\to\R$ is a formal psh function on $\cV$, 
  then $f\in\elde(g)$ iff
  there exists $\e>0$ such that $\log|f|\le(1+\e)g+A$ on $\cVqm$.
\end{Cor}
\begin{proof}
  This follows from the definition of $\elde(g)=\elde(Z(g))$
  and the strong continuity of $g$, $A$ and $\log|f|$ on $\cVqm$.
\end{proof}
As a second consequence of Theorem~\ref{thm:approx}
and~\ref{thm:subadditivity} 
we obtain the following approximation result.
\begin{Prop}\label{p:approx}
  For any formal psh function $g$ on $\cV$, there exists a countable
  sequence $(\fa_k)_1^\infty$ of primary ideals 
  and real numbers $c_k>0$ such
  that $g_k:=c_k\log|\fa_k|$ decreases to $g$ on all of $\cV$ as
  $k\to\infty$.
\end{Prop}
\begin{proof}
  Set $\fa_k=\elde(2^kg)+\fm^{4^k}$ and $c_k=2^{-k}$. By
  Theorems~\ref{thm:approx} and~\ref{thm:subadditivity}, 
  the sequence
  $h_k:=2^{-k}\log|\elde(2^kg)|$ decreases to $g$, hence
  $g_k=\max(h_k, -2^k)$ also decreases to $g$ on $\cVqm$. 
  We now show that this property
  automatically extends to $\cV\setminus\cVqm$. Let $\nu$ belong to
  the latter set, and pick $t>g(\nu)$ (which might be $-\infty$). 
  By the definition of $g(\nu)$, there exists $\pi$ 
  large enough so that
  $g(r_\pi(\nu))<t$, and then $g_k(r_\pi(\nu))<t$ for $k$ large
  enough. But we also have $g(\nu)\leq g_k(\nu)\leq g_k(r_\pi(\nu))$
  since $g_k$ is formal psh, and it follows that 
  $g_k(\nu)$ decreases to $g(\nu)$.
\end{proof}
We would like to stress the analogy with the usual 
complex case, where each psh function 
can be written as the decreasing limit of 
\emph{smooth} psh functions. 

Finally we show that formal psh functions satisfy the following
uniform Izumi-type bound.
\begin{Prop}\label{p:strict} For every
  $\nu\in\cVqm$, there exists $C=C(\nu)>0$ such that
  \[
  Cg(\nu_\fm)\leq g(\nu)\leq g(\nu_\fm)
  \]
  for each formal psh function
  $g$. In particular, either $\max_\cV g<0$ or $g\equiv 0$.
\end{Prop}
\begin{proof}
  The existence of $C$ and the left-hand inequality follow from Izumi's
  theorem~\cite{izumi,ELS}. 
  Indeed, the latter theorem yields a constant $C>0$ such that 
  $\nu(f)\le C\nu_{\fm}(f)$ for 
  all functions $f\in R$, hence $\log|\fa|(\nu)\ge C\log|\fa|(\nu_\fm)$ 
  for each ideal $\fa\subset R$, 
  and the result follows since every formal psh function 
  is the pointwise limit of functions of this type. The right-hand
  inequality 
  is the special case of $g\leq g\circ r_\pi$ when $\pi$ is the blow-up
  of $X$ at $0$: see Proposition~\ref{prop:fpshdecreasing}.
\end{proof}
%
%
\subsection{Nef envelopes}
The key technical tool for proving 
Theorem~\ref{thm:approx} and~\ref{thm:subadditivity} 
consists of approximating an arbitrary nef Weil divisor by
special ones, which are in turn well approximated from
below by Weil divisors associated to ideals.
\begin{Def}
  If $W$ is a Weil divisor on $\fX$, we let $\elin(W)$ be the
  ideal of functions $f\in\hat R$ such that $Z(f)\leq W$.
\end{Def}
Clearly $\elin(W)$ is an integrally closed ideal in $\hat R$, and by
definition $Z(\elin(W))\leq W$. In general, $\elin(W)$ may be
reduced to $0$. However, there is an important class of
Weil divisors for which $\elin(W)$ is nontrivial.
\begin{Def}
  If $\pi$ is any simple modification, a Weil divisor 
  $W$ on $\fX$ is said to be 
  \emph{determined on the dual complex $|\Delta(\pi)|$}
  if $W=Z(h)$ for some bounded function 
  $h:\cV\to\R$ such that $h = h \circ r_\pi$.
\end{Def}
For instance any Cartier divisor is determined on a dual
complex.

If $W=Z(h)$ is determined on $|\Delta(\pi)|$, then
$\elin(Z)\ne0$. Indeed, 
the condition $Z(f)\leq W$ only needs
to be tested on $|\Delta(\pi)|$ 
by Proposition~\ref{prop:fpshdecreasing} applied to $\log|f|$.
But $h$ is bounded so a sufficiently high power of $\fm$ is always
included in $\elin(W)$.

Since the cone of nef Weil divisors is closed and stable by taking max, 
the following definition makes sense.
\begin{Def}
  For any Weil divisor $W$ on $\fX$ 
  determined on some dual complex, we
  define the \emph{nef envelope} $\hat W$ as the supremum of all
  nef Weil divisor $Z$ on $\fX$ such that $Z\leq W$.
\end{Def}
If $f\in\hat R$ and $t>0$, $Z(f)\leq t W$ is equivalent to
$Z(f)\leq t\hat W$ since $t^{-1}Z(f)$ is nef. In other
words, we have $\elin(tW)=\elin(t\hat W)$ for all $t>0$.

The next result is in essence contained in~\cite{ELS}, see
Remark~\ref{rem:ELS} below.
\begin{Prop}\label{p:ELS}
  Let $W=\hat Z$ be the nef envelope of some Weil divisor $Z$
  determined on a dual complex.  For $k>0$ set
  $W_k:=\frac1k Z(\elin(kW))$. 
  Then there exists a constant $C>0$
  such that
  \begin{itemize}
  \item[(i)]
    $W\leq Z(\elde(W))\leq W+A$;
  \item[(ii)]
    $\elde(kW)\subset\elin((k-C)W)$ for all $k>C$;
  \item[(iii)]
    $W_k\leq W\leq(1-C/k)W_{k-C}$ for all $k>C$;
  \item[(iv)]
    $\elde(W)=\elde(W_k)$ for $k$ large enough.
  \end{itemize}
\end{Prop}
\begin{Remark}
  The proof will show that the constant $C$ in the statement can be
  taken to be $-\min_{|\Delta(\pi)|}A/g_W$ whenever $W$ is determined on
  $|\Delta(\pi)|$.  Note also that (i) is a special case of 
  Theorem~\ref{thm:approx} (which we shall prove in general in
  Section~\ref{sec:pf}.) 
\end{Remark}
\begin{Remark}\label{rem:ELS}
  Pick any simple modification $\pi$, and an effective divisor $D$ in
  $X_\pi$. 
  Then $\fa_t = \pi_*
  \cO_{X_\pi}(-tD)$ satisfies $\fa_t\cdot\fa_s\subset \fa_{t+s}$, hence 
  defines a graded family of ideals $\fa_\bullet$
  in the sense of~\cite{ELS}. 
  Set $W = \lim_t \frac1t Z(\fa_t)$. Then~(iv) above shows that the
  multiplier ideal of $\fa_\bullet$ as defined in op.cit.\ coincides with
  $\elde(W)$; whereas (iii) is exactly~\cite[Corollary~B]{ELS}.
\end{Remark}
\begin{Remark}\label{rem:fpshenv}
  In view of Proposition~\ref{p:extens} 
  we can similarly define the \emph{psh envelope} of any 
  function on $\cVqm$ determined on a dual complex. 
  Proposition~\ref{p:ELS} continues to hold in this context.
\end{Remark}
\begin{proof}[Proof of Proposition~\ref{p:ELS}]
  Throughout the proof we shall make use of the identification 
  of Weil divisors on $\fX$ with real-valued functions on $\cVdiv$.
  The result is trivial when $W=0$ so we shall assume
  $\nu_\fm(W)<0$. We start by proving
  \begin{Lemma}\label{lem:sup}
    With the notation above we have $W=\sup_k W_k$. 
  \end{Lemma}
  \begin{proof}
    For any $k$, $W_k$ is nef and $\le W$ by definition. For the
    reverse inequality, let $\pi\in\fB$ be large enough, so that $Z$
    is determined in $|\Delta(\pi)|$. Since $W$ is nef, given $\e>0$,
    there exists by definition a primary ideal $\fa$ of $R$ and an
    integer $k>0$ such that $W$ and $\frac1k Z(\fa)$ are $\e$-close
    on $|\Delta(\pi)|\cap\cVdiv$.  Since $W<0$, we can find an
    ideal of the form $\fb = \fa^N\fm^M$ and $l>0$ such that
    $W-\e\le\frac1l Z(\fb)\le W$ on $|\Delta(\pi)|\cap\cVdiv$. 
    A fortiori we
    have $\frac1l Z(\fb)\leq Z$ on $|\Delta(\pi)|$, hence everywhere
    on $\cV$ since $Z=r_\pi^*Z$. This means that
    $\fb\subset\elin(lZ)=\elin(lW)$, \ie 
    $\frac1l Z(\fb)\leq W_l\leq W$.  
    As $\frac1l Z(\fb)$ can be made arbitrarily close to $W$ on
    $|\Delta(\pi)|\cap\cVdiv$, 
    we conclude that $W=\sup W_k$ on $|\Delta(\pi)|\cap\cVdiv$ for
    every large enough $\pi$, hence on all of $\cVdiv$.
 \end{proof}
  We continue the proof of Proposition~\ref{p:ELS}.  The right-hand
  inequality in~(i) follows by the definition of $\elde(W)$. On the
  other hand, as follows from the definition, the primary ideal
  $\elin(kW)$ is contained in its usual multiplier ideal
  $\cJ(\elin(kW))$, which is equal to $\elde(kW_k)$ 
  by Proposition~\ref{p:usual}. 
  By the usual subadditivity theorem (Theorem~\ref{t:standard})
  we get
  $\elde(kW_k)\subset\elde(W_k)^k$.
  Moreover, $\elde(W_k)\subset\elde(W)$ since $W_k\leq W$. 
  We deduce from this chain of inequalities that
  $Z(\elin(kW))\leq k Z(\elde(W))$, \ie
  $W_k\leq Z(\elde(W))$, and (i) follows from Lemma~\ref{lem:sup} by
  letting $k\to\infty$.

  We now prove (ii).  Pick $\pi$ such
  that $Z=r_\pi^*Z$. If $f\in\elde(kW)$, then $Z(f)\leq kW+ A$
  by definition.  Now we have $A\leq -CW$ on $|\Delta(\pi)|$ for some
  $C>0$ (which only depends on $\pi$) since $W<0$. We
  infer that $Z(f)\leq (k-C)W\leq(k-C)Z$ on $|\Delta(\pi)|$, hence
  everywhere since $Z= r_\pi^*Z$. This means that
  $f\in\elin((k-C)Z)=\elin((k-C)W)$ as claimed.
  
  To get (iii), we use (ii). By applying (i) to $kW$
  instead of $W$, we get 
  $kW\leq Z(\elde(kW))\leq Z(\elin((k-C)W))$, and the result
  follows. Finally, to get (iv), first note that
  $\elde(W)=\elde(W+\frac CkW_C)$ for $k$ large enough. 
  Now $kW_k\leq (k-C)W_{k-C}+CW_C$, so 
  $W+\frac Ck W_C\leq W_k$ by (iii), which concludes the proof.
\end{proof}
%
%
\subsection{Proofs of Theorems~\ref{thm:approx} 
and~\ref{thm:subadditivity}.} \label{sec:pf}
Both proofs proceed by a reduction to the case of nef Weil 
divisors obtained as nef envelopes of divisors determined 
on a fixed dual complex.
\begin{proof}[Proof of Theorem~\ref{thm:approx}]
  Start with an arbitrary nef Weil divisor $W$ on $\fX$. 
  For every $\pi$,
  we can consider the incarnation $W_\pi\in\div(\fX)$ 
  of $W$ on $X_\pi$, viewed as a Cartier divisor on $\fX$.
  In general, $W_\pi$ may not be nef;
  let $\hat W_\pi$ be its nef envelope. 
  Since $W$ is nef, we have
  $W\le\hat W_\pi\le W_\pi$ by construction. 

  On the other hand, $W_\pi$ and $W$
  coincide on the vertices of $\Delta(\pi)$, and it follows 
  that the net $\hat W_\pi$ decreases pointwise to $W$ as 
  $\pi$ tends to infinity in the
  directed set of simple modifications.  
  We will rely on the following two results:
  \begin{Lemma}\label{l:monotone}
    Let $(Z_i)_{i\in I}$ be a decreasing net of nef Weil divisors
    on $\fX$, indexed by a directed set $I$, and 
    set $Z=\inf Z_i$.
    Then $\elde(Z) = \bigcap_i \elde((1+\e)Z_i)$ for
    every $\e>0$ small enough.
  \end{Lemma}
  \begin{Lemma}\label{l:modulo}
    Let $\fa_i$ be a decreasing net
    of ideals in $\hat R$, indexed by a
    directed set $I$, and set $\fa= \bigcap_I\fa_i$. Then
    $\inf_{i\in I}Z(\fa_i) = Z(\fa)$.
  \end{Lemma}
  We deduce from these lemmas that
  $Z(\elde(W))=\lim_{\pi}Z(\elde((1+\e)\hat W_\pi))$ for $\e>0$ small
  enough. We conclude by applying Proposition~\ref{p:ELS} (i) to
  $(1+\e)\hat W_\pi$ which by definition is the nef envelope of a
  Cartier divisor.
\end{proof}
\begin{proof}[Proof of Lemma~\ref{l:monotone}]
  We have $\elde(Z) = \elde ((1+\e)Z)$ for $\e\ll 1$, hence 
  $\elde(Z)\subset \bigcap_i \elde((1+\e)Z_i)$ 
  since $Z\leq Z_i$. Conversely,
  $f\in \bigcap_i \elde((1+\e)Z_i)$ implies $Z(f) \le (1+\e) Z_i+A$
  for all $i$ hence $f \in \elde(Z)$ since $Z_i$ converges to $Z$.
\end{proof}
\begin{proof}[Proof of Lemma~\ref{l:modulo}]
  For every
  ideal $\fb$ of $\hat R$, $Z(\fb)$ is the decreasing limit of
  $Z(\fb+\fm^k)$ as $k\to\infty$. Now, because $\hat R/\fm^k$ has
  finite length, it is Artinian, so that the decreasing sequence of
  ideals $\fa_i+\fm^k$ is stationary for $i\geq i(k)$. 
  We claim that the stationary value is actually $\fa+\fm^k$. 
  Granting this for a moment, we get
  \[
  Z(\fa)\leq\inf_{i\in
    I}Z(\fa_i)\leq Z(\fa_{i(k)}+\fm^k)=Z(\fa+\fm^k),
  \] 
  and so
  $Z(\fa)=\inf_{i\in I}Z(\fa_i)$ follows by letting
  $k\to\infty$. To prove the claim, let $p_l:\hat R\to \hat R/\fm^l$
  denote the projection. Each quotient map $\hat R/\fm^{l+1}\to\hat
  R/\fm^l$ induces a surjection from 
  $\bigcap_{i\in I}p_{l+1}(\fa_i)\to\bigcap_{i\in I}p_l(\fa_i)$,
  because both sides are reductions of the same $\fa_i$ for 
  $i\in I$ large enough. Because of
  this, given $f_k\in\bigcap_{i\in I}\fa_i+\fm^k$, there exists
  $f\in\hat R$ such that $f=f_k$ mod $\fm^k$ and 
  $f\in\bigcap_{i\in I}(\fa_i+\fm^l)$ for every $l\geq k$. 
  But now $f\in\fa_i$ mod $\fm^l$ for every large $l$ implies
  $f\in\fa_i$ because ideals of $\hat R$ are
  closed by Krull's lemma, and so we get 
  $f\in\bigcap_{i\in I}\fa_i=\fa$. 
  This proves the claim and the lemma.
\end{proof}
\begin{proof}[Proof of Theorem~\ref{thm:subadditivity}]
 By the same token, the proof boils down to 
 \begin{Lemma}
   Let $(\fa_i)_{i\in I}$ and $(\fb_i)_{i\in I}$ 
   be decreasing nets of ideals in $\hat R$, indexed
   by a directed set $I$. 
   Then $\bigcap_i(\fa_i\fb_i)=(\bigcap_i\fa_i )(\bigcap_i\fb_i)$.
 \end{Lemma}
 Again, we only have to check it modulo $\fm^l$ for all $l$ in which
 case the proof is obvious since the $\fa_i$ are stationary modulo $\fm^l$. 
\end{proof}
%
%
%
%
%

\section{Intersection theory and Monge-Amp{\`e}re operator}\label{sec:inter}
%
%
\subsection{Intersection of nef Weil divisors}\label{sec:inter1}
Fix a modification $\pi\in\fB$. We denote the intersection number of
any $n$-tuple of divisors $Z_1, \dots, Z_n \in \div(\pi)$ by $\langle
Z_1, \dots , Z_n \rangle \in \R$.  If $Z_1 = \sum_E r_E E$, then
$\langle Z_1, \dots , Z_n \rangle = \sum_E r_E\,\langle Z_2|_E, \dots
, Z_n|_E \rangle$. By expanding further in terms of the coefficients
of the other divisors, we see that the intersection product is
multilinear and continuous.  It is also symmetric in all variables and
satisfies the following important monotonicity property:
\begin{Lemma}
  Suppose $Z'_i,Z_i\in \div(\pi)$ are nef divisors
  satisfying $Z'_i \le Z_i$. Then
  \[
  \langle Z'_1,\dots, Z'_n\rangle \le \langle Z_1,\dots,Z_n\rangle.
  \]
\end{Lemma}
\begin{proof}
  By symmetry, the proof reduces to the case 
  $Z'_i = Z_i$ for $i\ge2$. But is then clear, 
  since $\langle Z_1,\dots,Z_n\rangle
  =\langle Z'_1,\dots,Z'_n\rangle+\langle(Z_1-Z'_1),Z'_2,\dots,Z'_n\rangle$
  and $Z_1-Z'_1$ is effective while $Z'_i$ is nef for $i\geq 2$.
\end{proof}
We now extend the previous construction to $\Cdiv(\fX)$. For any
collection $Z_1,\dots,Z_n$ of Cartier divisors on $\fX$
determined on some $X_\pi$, we set
\begin{equation}\label{e:cartint}
  \langle Z_1, \dots , Z_n \rangle \= 
  \langle Z_{1,\pi}, \dots , Z_{n,\pi} \rangle.
\end{equation}
This definition does not depend on $\pi$, and we have
\begin{Prop}
  The intersection product $\Cdiv(\fX)^n\to\R$ defined 
  by~\eqref{e:cartint} is symmetric, multilinear and continuous 
  for the strong topology.  
  Moreover, if $Z'_i\le Z_i$ are nef Cartier divisors, then
  \begin{equation}\label{e:incr}
    \langle Z'_1,\dots,Z'_n\rangle \le \langle Z_1,\dots,Z_n\rangle.
  \end{equation}
\end{Prop}
Note that $\langle Z_1,\dots,Z_n\rangle\leq 0$ when the $Z_i$ are nef,
because $Z_i$ is then $\leq 0$.

Define the \emph{mixed multiplicity} 
$e\langle\fa_1,\dots,\fa_n\rangle$ of a
collection of $n$ primary ideals $\fa_i$, as
$e\langle\fa_1,\dots,\fa_n\rangle\=\dim R / (f_1,\dots,f_n)$ 
where $f_i$ is a generic element in $\fa_i$, see also~\cite{rees}. 
Then the following holds, see~\cite{rama}:
\begin{equation}\label{e:ideal}
  -\langle Z(\fa_1),\dots,Z(\fa_n) \rangle 
  = e\langle\fa_1,\dots,\fa_n\rangle.
\end{equation}

The inequality~\eqref{e:incr} allows us to extend the definition to
nef Weil divisors.
\begin{Def}
  Let $Z_1, \dots , Z_n$ be nef Weil divisors on $\fX$. 
  Then we set
  \[
  \langle Z_1, \dots , Z_n \rangle
  \= \inf\langle W_1, \dots , W_n
  \rangle \in [-\infty,0],
  \] 
  where the infimum is taken over
  nef Cartier divisors such that $W_i \ge Z_i$ for all $i$.
\end{Def}
By monotonicity, this definition is consistent with the intersection
defined above for nef Cartier divisors.
\begin{Prop}{~}\label{p:intersect}
  The intersection product on nef Weil divisors
  is symmetric, upper semicontinuous, and
  continuous along decreasing families.
  It is also $1$-homogeneous, additive,
  and increasing in each variable. 
  Further, $\langle Z_1,\dots,Z_n\rangle<0$ unless
  $Z_i=0$ for some $i$.
\end{Prop}
\begin{Remark}\label{rem:essential}
  Our approximation result Proposition~\ref{p:approx} whose proof is
  based on multiplier ideals is used here in an essential way to prove
  the additivity of the intersection product.
\end{Remark}
\begin{proof}
  Symmetry, homogeneity and monotonicity are all clear by
  definition. Suppose that all the $Z_i$ are non-zero. 
  Then there exists
  $\e>0$ such that $Z_i\le\e Z(\fm)$ for all $i$ by
  Proposition~\ref{p:decrease}, and it follows that 
  $\langle Z_1,\dots,Z_n\rangle\leq -\e^n$. 

  Let us prove upper semicontinuity. 
  The proof relies on the following remark: if $W$ is
  a Cartier divisor on $\fX$, 
  then the set of nef Weil divisors $Z$ on $\fX$ such that
  $Z<W$ (coefficient-wise) is open in the nef cone in 
  $\div(\fX)$. Indeed, this
  condition needs only be tested on the finitely many exceptional
  primes of a determination $\pi$ of $W$ by
  Proposition~\ref{p:decrease}. 
  Now let $Z_1,\dots,Z_n$ be given nef Weil
  divisors. Since the intersection product of nef Weil divisors is
  always nonpositive, upper semicontinuity at $Z_1,\dots,Z_n$ is
  automatic if one the $Z_i$ is zero. We can 
  therefore assume that they are all non-zero, so that 
  $\langle Z_1,\dots,Z_n\rangle<0$ by what precedes.

  Let $t<0$ be a given real number such that
  $\langle Z_1,\dots,Z_n\rangle<t$. By definition, there
  exist nef Cartier divisors $W_i\geq Z_i$ such that 
  $\langle W_1,\dots,W_n\rangle<t$. 
  The $W_i$ are in particular non-zero, thus $W_i<0$
  coefficient-wise by Proposition~\ref{p:decrease}. 
  Therefore, upon
  replacing $W_i$ by $(1-\e)W_i$ for $\e\ll 1$ we can assume that
  $Z_i<W_i$ by the remark above.
  Thus the set $U_i$ of nef Weil divisors
  $Z'_i$ such that $Z'_i<W_i$ is a neighborhood of $Z_i$ in the nef
  cone, and 
  $\langle Z'_1,\dots,Z'_n\rangle\leq\langle W_1,\dots,W_n\rangle<t$ 
  for $Z'_i$ in $U_i$. 
  This proves upper semi-continuity. Since the
  intersection product is increasing in each variable, continuity along
  decreasing families follows immediately.

  As to additivity, let $Z'_1$ be another nef
  Weil divisor. By Proposition~\ref{p:approx}, 
  there exists decreasing
  sequences $(W_{1,k})_{k=1}^\infty$, 
  $(W'_{1,k})_{k=1}^\infty$ 
  and $(W_{i,k})_{k=1}^\infty$, $i\geq 2$ of nef
  Cartier divisors on $\fX$, 
  converging to $Z_1$, $Z'_1$ and $Z_i$
  respectively. Thus additivity follows from the continuity 
  just proved 
  and from the additivity of the intersection pairing on 
  Cartier divisors.
\end{proof}
%
%
\subsection{Monge-Amp{\`e}re measures on $\cV$}\label{sec:inter2}
We now interpret the intersection product on nef Weil divisors
constructed above as a positive measure on $\cV$. 

Recall that a \emph{Radon measure} on $\cV$ is a 
continuous linear form on
the space $C^0(\cV)$ of continuous real-valued functions on $\cV$
(with respect to its weak topology), endowed with the $\sup$ norm. If
$W$ is a dense vector subspace of $C^0(\cV)$ that contains the
constant functions, then any non-negative linear form $L$ on $W$
defines a positive Radon measure on $\cV$, of total mass $L({\bf 1})$
with ${\bf 1}$ the constant function on $\cV$. Indeed, for every $f\in
W$, we have $(\inf f){\bf 1}\leq f\leq(\sup f){\bf 1}$, and thus
$(\inf f) L({\bf 1})\leq L(f)\leq(\sup f) L({\bf 1})$. In other words
$|L(f)|\leq L({\bf 1})\sup |f|$, and the result follows by uniform
continuity.
\begin{Prop}
  Suppose $g_1,\dots,g_{n-1}$ are formal psh functions on $\cV$,
  with associated nef divisors $Z_i$ on $\fX$.  
  If $\langle Z_1,\dots,Z_{n-1},Z(\fm)\rangle>-\infty$,
  there exists a unique positive Radon measure
  $\MA(g_1,\dots,g_{n-1})$
  on $\cV$ such that
  \[
  \int_\cV g\MA(g_1,\dots,g_{n-1})
  =\langle Z(g),Z_1,\dots,Z_{n-1}\rangle,
  \]
  for any formal psh function $g$ on $\cV$.
\end{Prop}
When $\langle Z_1,\dots,Z_{n-1},Z(\fm)\rangle>-\infty$ holds, 
we will say that
the \emph{Monge-Amp{\`e}re measure} $\MA(g_1,\dots,g_{n-1})$ is
well-defined. This is the case if the $g_i$ are bounded on
$\cV$. Indeed, $g_i\geq -C$ means $Z_i\geq CZ(\fm)$, thus $\langle
Z_1,\dots,Z_{n-1},Z(\fm)\rangle\geq -C^{n-1}$ by monotonicity.
\begin{proof} 
  If $W$ is a given nef Cartier 
  divisor on $\fX$, with associated continuous 
  (and piecewise affine) function $g_W$ on $\cV$,
  then $g_W$ is bounded on the compact space $\cV$, 
  \ie $g_W\geq -C$ for some $C>0$.  
  In terms of divisors, this means that $W\geq CZ(\fm)$,
  and thus our assumption implies that 
  $\langle Z_1,\dots,Z_{n-1},W\rangle>-\infty$ for any 
  nef Cartier divisor $W$. 
  Now if $Z$ is an arbitrary Cartier divisor on $\fX$, 
  we can write it as a
  difference of two nef Cartier divisors $Z=W_1-W_2$, and we then set
  $\langle Z_1,\dots,Z_{n-1},Z\rangle:=\langle
  Z_1,\dots,Z_{n-1},W_1\rangle-\langle
  Z_1,\dots,Z_{n-1},W_2\rangle$. This does not depend on the choice of
  $W_1$ and $W_2$ by additivity. We thus get a non-negative functional
  $Z\mapsto\langle Z_1,\dots,Z_{n-1},Z\rangle$ on $\Cdiv(\fX)$. By
  Proposition~\ref{p:cartier}, $\Cdiv(\fX)$ can be identified with the
  dense subspace $\PA(\cV)$ of $C^0(\cV)$, 
  see Proposition~\ref{p:dense}. 
  By the preceding discussion, this
  concludes the proof.
\end{proof}
By translating Proposition~\ref{p:intersect},~\eqref{e:incr}
and~\eqref{e:ideal} into statements about psh functions we obtain 
the following two results.
\begin{Prop}
  For any primary ideals $\fa_1,\dots,\fa_n$ of $R$, we have
  \[
  e\langle\fa_1,\dots,\fa_n\rangle=\int_\cV
  -\log|\fa_1|\,\MA(\log|\fa_2|,\dots,\log|\fa_n|).
  \]
\end{Prop}
\begin{Prop}
  Suppose $(g_{i,k})_{k=1}^\infty$ 
  is a sequence of formal psh functions
  decreasing to a formal psh function $g_i$ for $i=1,\dots,n-1$. 
  If $\rho\=\MA(g_1,\dots,g_{n-1})$ is defined, then so is
  $\rho_k\=\MA(g_{1,k},\dots,g_{n-1,k})$ for all $k$, 
  and $\rho_k\to\rho$ weakly as $k\to\infty$
\end{Prop}
Let us spell out more explicitly what happens for 
piecewise affine functions. 
\begin{Prop}\label{p:suppma}
  If $g_1,\dots,g_{n-1}$ are piecewise affine 
  formal psh functions on $\cV$,
  then $\MA(g_1,\dots,g_{n-1})$ is a finite sum of Dirac masses at
  divisorial valuations.
\end{Prop} 
\begin{proof} 
  Choose $\pi$ such that each $g_i$ belongs to $\Aff(\pi)$, 
  and let $g$
  be a piecewise affine test function, lying in $\Aff(\pi')$ for some
  $\pi'$ possibly very large. By definition, we then have
  $\int_{\cV}g\,\MA(g_1,\dots,g_{n-1})=\langle
  Z(g),Z(g_1),\dots,Z(g_{n-1})\rangle$, which can be computed 
  on $X_\pi$
  itself by the projection formula, even if $Z(g)$ is only 
  determined by $\pi'$. 
  The integral is thus equal to 
  $\sum_Eg(\nu_E)b_E\langle E,Z(g_1),\dots,Z(g_{n-1})\rangle$, 
  where the sum is over the vertices of
  $\Delta(\pi)$, $\nu_E=b_E^{-1}\ord_E$ is the normalized valuation
  proportional to $\ord_E$, and $b_E=\ord_E(\fm)$ is the normalization
  factor. We thus see that $\MA(g_1,\dots,g_{n-1})$ 
  is supported on the
  vertices of $\Delta(\pi)$, with mass 
  $b_E\langle E,Z(g_1),\dots,Z(g_{n-1})\rangle$ at 
  the vertex $\nu_E$ determined by
  a $\pi$-exceptional prime $E$.
\end{proof} 

In fact, this reasoning shows that the intersection number 
$\langle Z_1,\dots,Z_n\rangle$ of nef Weil divisors on $\fX$ 
is always finite when all
but at most one $Z_i$ are Cartier. Equivalently, the Monge-Amp{\`e}re
measure $\MA(g_1,\dots,g_{n-1})$ of formal psh functions is always
well-defined when all but at most one of the $g_i$ are piecewise affine.

In the case $n=2$, we thus see that $\MA(g)$ is always defined when $g$ is
a formal psh function on the valuative tree $\cV$. One can check that 
$\MA(g)$ coincides with the Laplacian $\Delta(-g)$ of the 
(positive) tree potential $-g$ as defined in~\cite{valtree}. 
In this special case, the Laplace operator
can be constructed purely in terms of the metric tree structure of
$\cV$. See also~\cite{BR1,thuillier1}. 
In the higher dimensional case, the metric tree structure is
replaced by the piecewise integral structure on $\cV$ introduced in
Section~\ref{sec:affine}, but this is \emph{not} enough to define the
Monge-Amp{\`e}re operator.
%
%
%
%
%
\section{Singularities of plurisubharmonic functions}\label{sec:psh}
Throughout this section, we work over the field of complex numbers,
and denote by $\cO$ the ring of germs of holomorphic functions at
the origin in $\C^n$. 
It lies between the two rings $R\subset\hat R$ we have
been considering, and primary ideals of all three rings can (and
will) be identified. Given an analytic ideal $\fa\subset\cO$ and a
valuation $\nu\in\hcV$, we will compute $\nu(\fa)$ as the value of
$\nu$ on the formal ideal $\fa\hat R$.
%
%
\subsection{Definitions}
Let $u$ be a \emph{psh germ}, that is, a plurisubharmonic function
defined in a neighborhood of the origin in $\C^n$. The origin is a
\emph{singularity} of $u$ if $u(0)=-\infty$.  When $u$ and $v$ are two
psh germs, we say that $u$ is \emph{more singular} than $v$ 
if $u\leq v+O(1)$ near $0$. 
We say that $u$ has \emph{singularities described by
$\fa^c$} for some ideal $\fa$ of $\cO$ and $c>0$ if $u=c\, \max \log
|f_k|+O(1)$ for some (hence any) choice of local generators $(f_k)$ of
$\fa$.

Recall that the \emph{Lelong number} $\nuL(u,0)\in \R_+$ of $u$ is the
supremum of all $c\geq 0$ such that $u$ is more singular than
$c\log|z|$. The supremum is actually attained, as follows from the
basic fact that $\sup_{\{\log|z|<t\}} u$ is a convex function of
$t$. When $u$ has singularities described by an analytic ideal $\fa$,
$\nuL(u,0)$ is just $\nu_\fm(\fa)$.

The \emph{multiplier ideal} $\cJ(u)\subset\cO$ of a psh germ $u$ is
the ideal of germs of holomorphic functions $f$ such that $fe^{-u}$ is
$L^2$ near $0$. When the singularity of $u$
is described by $\fa^c$, it is elementary to see that
this definition coincides with the one
given in Section~\ref{sec:def-mult}: 
see \eg~\cite[Proposition~1.7]{DK}.

For any psh germ, the multiplier ideals $\cJ(ku)$ allow us to
approximate the singularity of $u$ very precisely as $k\to\infty$. We
refer to~\cite[Theorem~4.2]{DK} for details on the following construction.
Fix a ball $B$ centered at the origin
on which $u$ is defined, and let $\cH (u,B)$ be the set of holomorphic
functions $f$ on $B$ such that 
$\| f \|^2 \= \int_B |f|^2 e^{-2u} < + \infty$. 
This is a Hilbert space, and one can show that the germs at $0$
of the functions of $\cH(u,B)$ generate $\cJ(u)$. This implies that
the psh germ 
$\Reg{u} := \sup \{\, \log |f|, \, f \in \cH(u,B),\, \|f\| =1\}$ 
on $B$ has singularities described by $\cJ(u)$ at $0$.

Now it is a fundamental fact that there exists a constant $C>0$
depending only on the dimension $n$ and the diameter of $B$,
such that
\begin{equation}\label{e:OT}
  u(p) - C \le \Reg{u}(p) \le \sup_{B(p,r)} u + C - n \log r,
\end{equation}
for all $p$ and all $r$ such that $B(p,r)\subset B$. 
Here the right hand
side is a consequence of the mean value inequality for $|f|^2$ when
$f$ is holomorphic, whereas the left hand-side is a deep fact, an
application of the Ohsawa-Takegoshi theorem.

Using~\eqref{e:OT}, it is easy to see that $u_k\=\frac1k\Reg(ku)$
converges pointwise to $u$ as $k\to\infty$, and that
$\frac1k\nu_\fm(\cJ(ku))=\nuL(u_k,0)$ converges to $\nuL(u,0)$. We
will now show that this convergence extends to $\cV$.
%
%
\subsection{Valuative transform of a psh germ}\label{sec:trace}
We can associate to a psh germ $u$ a Weil divisor $Z(u)\in\div(\fX)$
as follows. If $E$ is an exceptional prime on some $X_\pi$, we
define $-\ord_E Z(u)$ to be the generic Lelong number of $u\circ\pi$
along $E$, \ie its Lelong number at the very general point of
$E$. This only depends on the divisorial valuation $\ord_E$.
Equivalently, $[Z(u)_\pi]$ is the divisorial part in the Siu 
decomposition of the current $dd^c(u\circ\pi)$ on $X_\pi$.
\begin{Thm}\label{thm:psh-nef}
  For any psh germ $u$, the Weil divisor $Z(u)$ on $\fX$ is nef.
\end{Thm}
\begin{Def}
  We write $\hat u$ for the associated formal psh function on $\cV$
  and call it the \emph{valuative transform} of $u$.  
\end{Def}
Let us relate the analytic multiplier ideal $\cJ(u)$ 
of $u$ to its valuative counterpart $\elde(\hat u)$.
For that purpose, note that the family of ideals $\cJ((1+\e)u)$
increases as $\e>0$ decreases to $0$, thus has to be stationary by
the Noetherian property.  This justifies the following definition.
\begin{Def}
  For any psh germ $u$, we set $\cJ_+(u)\= \cJ((1+\e)u)$ for $\e>0$
  sufficiently small.
\end{Def}
\begin{Remark}\label{rem:open}
  The Openness Conjecture by
  Demailly and Koll{\'a}r~\cite[Remark~5.3]{DK}
  asserts that $\cJ_+(u)=\cJ(u)$ for any psh germ $u$.
  In~\cite{FJ-multsing} this conjecture was settled affirmatively
  in dimension two.
\end{Remark}
\begin{Thm}\label{thm:mult-comp}
  If $u$ is any psh germ, then we have
  $\elde(\hat u)\cap\cO=\cJ_+(u)$.
\end{Thm}
Let us give an explicit interpretation of the restriction of
the valuative transform $\hu$ to a simplex in $\cV$.
Let $\pi$ be a simple modification, and $E_1,\dots,E_k$ be
exceptional primes of $\pi$ determining a simplex $|\sigma|$ in 
$|\Delta(\pi)|\subset\cV$. Let $z_1,\dots,z_k$ be 
local equations of $E_1,\dots,E_k$ at a generic point 
$p$ on $E_1\cap\dots\cap E_k$, and complete these to 
local coordinates $z_1,\dots,z_n$ at $p$. 
Then $|\sigma|$ consists of
those valuations $\nu_{z,w}$, $w\in\R_+^n$ that are monomial in this
coordinate system, and such that $\sum_{1\le i\le k} b_i w_i=1$,
$w_i=0$ for $i>k$. Here $b_i=\ord_{E_i}(\fm)$. We will
identify $\sigma$ with the set of such $w\in\R_+^n$.

Now for any $w\in\R_+^n$, if $v$ is a psh function on $X_\pi$ near
$p$, set 
\begin{equation*}
  \nuK_{z,w}(v) 
  = \sup \{ t\ge0,\, v \le t\, \log \max |z_i|^{w_i} + O(1)\}.
\end{equation*}
This weighted version of the Lelong number was
introduced in~\cite{kis} and is referred to as the \emph{Kiselman
number}.  Kiselman numbers are $1$-homogeneous and monotonous in $w$:
$\nuK_{z,c w} = c \nuK_{z,w}$ for any $c>0$; and $w_i\le w_i'$ for
all $i$ implies $\nuK_{z,w}(v)\ge \nuK_{z,w'}(v)$.  In particular, the
function $w \mapsto \nuK_{z,w}(v)$ is continuous in $w$.

When $w\in \Q_+^n$, the monomial valuation $\nu_{z,w}\in|\sigma|$ 
is divisorial
so we may write $\nu_{z,w}= c\ord_E$ for some exceptional divisor
$E$ of a modification $\mu$ of $X_\pi$ above $\pi^{-1}(0)$, 
and some positive real number $c>0$.  
It is then not difficult to see that 
$\nuK_{z,w}(v) = c\ord_E(v \circ \mu)$. 
If $v = u \circ \pi$, this shows
$\hu(\nu_{z,w})=-\nuK_{z,w}(u \circ \pi)$ when $w\in \Q_+^n$.  By
continuity, we conclude
\begin{Prop}
  For any simple modification $\pi$, and any 
  quasimonomial valuation $\nu_{z,w}$ in a simplex 
  $|\sigma|\subset\cV$ appearing as a face on $|\Delta(\pi)|$, 
  and any psh germ $u$, the number $\hu(\nu_{z,w})$ equals the
  Kiselman number $-\nuK_{z,w}(u \circ \pi)$. 
\end{Prop}
In particular, the Kiselman numbers are concave functions of
the weight $w$, see Proposition~\ref{p:extens}. 
This fact is due to Kiselman.
From~\cite[Proposition~3.12]{dem2}, one also infers 
\begin{Cor}
Let $u_j, u$ be psh germs such that $u_j \to u$ in $L^1_\mathrm{loc}$.
Then
\begin{equation}\label{e:semi}
\liminf_{j\to\infty} \hat{u}_j(\nu)  \ge \hat{u}(\nu)
\end{equation}
for any quasimonomial valuation $\nu$.
\end{Cor}

\begin{proof}[Proof of Theorem~\ref{thm:psh-nef}]
  By what was recalled above, 
  $u$ is more singular than the psh germ $\Reg(u)$, 
  which has a singularity described by $\cJ(u)$.
  Hence $Z(u)\leq Z(\cJ(u))$. On the other hand, let $f\in \cJ(u)$. 
  If $\pi\in\cB$ and
  $E$ is an exceptional prime of $\pi$, the local integrability of
  $|f|^2 e^{-2u}$ implies the integrability of $|f\circ \pi|^2
  e^{-2u\circ \pi} |J\pi|^2$ on $X_\pi$ by the change of variables
  formula, and this in turn is easily seen to imply 
  $\ord_E(f \circ \pi)+ A(\ord_E) >\ord_E(u \circ \pi)$, 
  as a consequence of Fubini's
  formula at a generic point of $E$. It follows that 
  $\ord_E(\cJ(u)) + A(\ord_E) \ge \ord_E(u)$, and we conclude
  \begin{equation}\label{e:app}
    Z(\cJ(u)) \ge Z(u) \ge Z(\cJ(u)) - A.
  \end{equation}
  It follows that $Z(u)$ is the limit of the nef divisors
  $\frac1kZ(\cJ(ku))$, hence is nef.
\end{proof}
\begin{proof}[Proof of Theorem~\ref{thm:mult-comp}]
  If $f\in\cJ_+(u)$, then we have $Z(f)\leq Z((1+\e)u)+A$
  for $\e>0$ small enough by~\eqref{e:app}, and this 
  implies that
  $f\in\elde(\hat u)$ by definition.
  
  In case $u$ has singularities described by $\fa^c$ for some analytic
  ideal $\fa$ of $\cO$ (not necessarily primary), 
  the reverse inclusion is proved exactly as in
  Proposition~\ref{p:usual}, using an analytic log-resolution
  of $\fa$.
  
  Now let $u$ be an arbitrary psh germ, defined on a small ball
  $B$ centered at $0$, and let $u_k=\frac1k\Reg(ku)$, so that 
  $u_k$ has
  singularities described by $\cJ(ku)^{1/k}$. We have
  $\elde(\hu)=\elde((1+\e)\hu)$ for small $\e>0$ by definition, 
  and $\elde((1+\e)\hu)\subset\elde((1+\e)\hu_k)$ for all $k$ 
  since $u_k$ is less singular than $u$. 
  By the special case just considered, we
  have $\elde((1+\e)\hu_k)\cap\cO=\cJ((1+\e)u_k)$. Thus if
  $f\in\elde(\hu)\cap\cO$, we have $\int |f|^2e^{-2(1+\e)u_k}<+\infty$
  for all $k$. On the other hand, $u_k$ satisfies $\int
  e^{2k(u_k-u)}<+\infty$, and these two inequalities imply
  $\int|f|^2e^{-2u}<+\infty$ by the H{\"o}lder inequality, as
  in~\cite[Proof of Theorem~4.2]{DK}. This concludes the proof.
\end{proof}
%
%
\subsection{Weights}
By a \emph{psh weight}, we will mean a psh germ $\varphi$ having
an isolated singularity at the origin and
such that $e^\varphi$ is continuous.

We now introduce a class of psh germs that will play an
important role in the sequel.
\begin{Def}\label{def:reg} 
  A psh germ $\varphi$ is said to be 
  \emph{tame} if there exists a constant $C>0$ such that for every
  $t>0$ and every germ $f\in\cJ(t\varphi)$, $\log|f|$ is more singular
  than $(t-C)\varphi$.
\end{Def}
The definition can be paraphrased by saying that the $L^2$ and
$L^{\infty}$ conditions for holomorphic functions with respect to
$e^{-t\varphi}$ stay uniformly bounded as $t\to\infty$. Concretely, we
have the following lemma.
\begin{Lemma}\label{lem:easier1}
  Suppose $\varphi$ is a tame psh germ, and $\varphi_k$ is any psh
  function with singularities described by $\frac1k \cJ(k\varphi)$.
  Then $\varphi \le \varphi_k \le (1-C/k)\varphi$ up to bounded
  functions.
\end{Lemma}
\begin{proof}
  It is sufficient to prove the lemma for $\varphi_k \= \frac1k \Reg
  (k\varphi)$. From~\eqref{e:OT}, 
  we get $\varphi\le\varphi_k+O(1)$. 
  Conversely, pick a finite set of generators $f_i$ of the
  ideal $\cJ(k\varphi)$. Then 
  $\varphi_k =\frac1k \log \max|f_i|+O(1)$ 
  near the origin and from the definition of tame
  weight we have $\log |f_i| \le (k-C) \varphi$ up to a constant.
\end{proof}
We shall be mostly concerned with tame weights.
Here is an important class of such weights.
\begin{Lemma} 
  Let $\varphi$ be a psh germ such that $e^\varphi$ is 
  $\alpha$-H{\"o}lder for some $\alpha>0$.  Then $\varphi$ is a tame weight. In
  fact, we can take $C=n/\alpha$ in the above definition.
\end{Lemma} 
\begin{proof} 
  Assume that $\varphi$ is defined on a ball $B$. 
  By~\eqref{e:OT} we have
  $\frac1k\Reg(k\varphi)(p)\leq\sup_{B(p,r)}\varphi-n\log r+C$ 
  for every
  $r$ small enough and some uniform constant $C>0$. Now $e^\varphi$ is
  $\alpha$-H{\"o}lder, thus $\sup_{B(p,r)}e^\varphi\leq
  e^{\varphi(p)}+Cr^{\alpha}$. Choosing $r=e^{\varphi(p)/\alpha}$ now
  shows that $\sup_{B(p,r)}\varphi\leq \varphi(p)+O(1)$, and this
  concludes the proof.
\end{proof}       
\begin{Remark}\label{rmk:tame}
  One can export the notion of tameness into $\cV$.  
  Define a formal psh function $g$ on $\cV$ to be \emph{tame} when
  $\elde(kg)\subset\elin((k-C)g)$ for some $C>0$ and all $k$, 
  or equivalently $g \le h_k \le (1-C/k)g$ 
  with $h_k:=\frac1k\log|\elde(kg)|$. 
  In particular, $\elin(kg)$ is non-trivial, and in fact 
  $g_k:=\frac 1k\log|\elin(kg)|$ converges to $g$ as $k\to\infty$. 
  Proposition~\ref{p:ELS} and Remark~\ref{rem:fpshenv} 
  assert that the
  formal psh envelope of a function determined on a dual 
  complex is tame, and the rest of Proposition~\ref{p:ELS}
  in fact holds for any tame formal psh function.
\end{Remark}
%
%
\subsection{Relative types and valuations}
Let us say that a psh weight $\varphi$ is \emph{maximal} if it is
maximal outside the origin, that is, $(dd^c\varphi)^n=\lambda\delta_0$
for some $\lambda>0$.  If $u$ is an arbitrary psh germ,
Rashkovskii~\cite{rash} defined the \emph{relative type}
$\sigma(u,\varphi)$ of $u$ with respect to $\varphi$ as the largest
$c\geq 0$ such that $u$ is more singular than $c\varphi$.  It is a
fundamental fact that $\sup_{\{\varphi\leq t\}}u$ is a convex function
of $t$ under the maximality assumption.
Hence $\sigma(u,\varphi)<+\infty$ and that $u$
is more singular than $\sigma(u,\varphi)\varphi$.

We now relate quasi-monomial valuations to relative types with respect
to tame maximal weights.
\begin{Thm}\label{thm:relat}
  If $\varphi$ is a tame maximal weight and $u$ is any psh germ, 
  then $u\le   \varphi+O(1)$ iff $\hat u\le \hat\varphi$. 
  In other words for   every psh germ $u$, we have
  $\sigma(u,\varphi)=\sup \{ c\ge0,\, \hat u \le c\,\hat\varphi\}$.
\end{Thm}
\begin{Thm}\label{thm:repr} 
  For any quasi-monomial valuation $\nu\in\cV$, 
  there exists a tame maximal psh
  weight $\varphi$ such that $\hat{u}(\nu) = \sigma(u,\varphi)$
for any psh function $u$.
 
The weight $\varphi$ is unique up to a bounded
  function, and $\hat\varphi$ can be characterized as 
  the largest formal psh 
  function $g$ on $\cV$ such that $g(\nu)=-1$.
\end{Thm}
\begin{Remark}
  The proof of Theorem~\ref{thm:repr} shows that the constant $C$ making
  $\varphi$ tame can be taken to be $A(\nu)$, 
  and one can show that this is the best
  possible choice. Also a substantial part of the result is a direct
  consequence of~\cite[Theorem~4.1]{rash}. Our main contribution lies in
  proving that $\varphi$ is continuous and tame.
\end{Remark}
\begin{proof}[Proof of Theorem~\ref{thm:relat}]
  If $u$ is more singular than $\varphi$, 
  then of course $\hat u\leq\hat\varphi$. 
  To get the converse, let $\varphi_k$ be a psh germ with singularities
  described by $\cJ(k\varphi)^{1/k}$. 
  If $\hu\le\hvarphi$, 
  then a fortiori $\hu\le\frac1k\log|\cJ(k\varphi)|$. 
  Pick a log-resolution $\pi_k$ of the primary ideal $\elde(k\varphi)$, 
  so that $\varphi_k\circ\pi_k$ has divisorial singularities. 
  Then we have 
  \[
  \ord_E(u\circ\pi_k)
  =-\hu(\ord_E)\geq-\frac1k\log|\cJ(k\varphi)|(\ord_E)
  =\ord_E(\varphi_k\circ\pi_k)
  \]
  for each exceptional prime $E$ of $\pi_k$, and it follows from 
  a well-known   theorem of Siu that $u\circ\pi_k-\varphi_k\circ\pi_k$ 
  is psh up to a bounded term. It is in particular bounded from above, 
  and we infer that $u\le\varphi_k+O(1)$ near the origin. 
  Now Lemma~\ref{lem:easier1} implies that $u\le (1-C/k)\varphi+O(1)$, 
  \ie $\sigma(u,\varphi)\geq 1-C/k$ for all $k$, and we 
  conclude that $\sigma(u,\varphi)\ge 1$, \ie $u\le\varphi+O(1)$. 
\end{proof}
\begin{proof}[Proof of Theorem~\ref{thm:repr}]
We may suppose we work on the unit ball $B$, and define 
$\varphi$ using a classical envelope construction: 
$$\varphi = \sup\{ u\le 0 \text{ psh on } B \text{ such that } \hat{u}(\nu) \le -1 \}.$$ 
Choquet's lemma gives a sequence of negative psh functions $v_j$ with $\hat{v}_j(\nu) \le -1$ and 
increasing to a function $v$ whose upper-semicontinuous regularization $v^*$ is a psh function equal to $\varphi^*\ge \varphi$.  From~\eqref{e:semi}, we infer $\hat{v}^*(\nu) \le \liminf \hat{v}_j(\nu) \le -1$.
This shows $v^* \le \varphi$, hence $\varphi= \varphi^*$ is psh.

As $\log |z|$ is a negative psh function whose valuative transform is constant equal to $-1$, we get
$\log|z| \le \varphi$ and  $\hat{\varphi}(\nu) =-1$.
It is not difficult to see that $\log |z|$  is also the largest negative psh function on $B$ 
such that $\hat{u}(\nu_\fm) \le -1$. Whence $\log |z| \le \varphi \le c \log |z|$ with $c = -\hat{\varphi}(\nu_\fm)$. This implies $\varphi$ to be locally bounded outside $0$ and $e^\varphi$ to be continuous at $0$. Maximality of $\varphi$ follows from standard arguments.

We now show that 
\begin{equation}\label{e:form}
\hat{\varphi} = \sup \{ h \text{ formal psh}, \, h(\nu) \le -1 \}~.
\end{equation}
Denote by $g$ the right hand side of~\eqref{e:form}.
As $\hat{\varphi}$ is a formal psh function with $\hat{\varphi}(\nu) = -1$, we have $\hat{\varphi} \le g$.
Set $g_k = \frac{1}{k} \elin(kg)$. We claim that $g_k \to g$. To see that, pick a simple blow-up $\pi$  such that $\nu \in | \Delta (\pi)|$. Define the bounded function $h: \cV \to \R$ by setting 
$h(\nu) = -1$, $h \equiv 0$ on $|\Delta(\pi)| \setminus \{ \nu \}$ and $h = h \circ r_\pi$. Then $g$ is the formal psh function associated to the nef envelope of $Z(h)$ so that we may apply Proposition~\ref{p:ELS}. We conclude $g_k \to g$.

Now the ideals $\elin(kg)$ are primary. One can hence find a negative psh function $u_k$ on $B$ with singularities described by $\frac{1}{k}\elin(kg)$. As $g_k \to g$, we get $(1+\e) u_k \le \varphi$ for large $k$ and $\e$ arbitrarily small. Then $(1+\e)g_k \le \hat{\varphi}$ and thus $ g \le \hat{\varphi}$.
This proves~\eqref{e:form}.

\smallskip

Next we show the continuity of $e^\varphi$. The argument is standard. 
By what precedes, we only need to check the continuity of $\varphi$ outside $0$. Write $\varphi_k = \frac{1}{k} \Reg (k\varphi)$. We shall prove that $\varphi_k \to \varphi$ uniformly on any compact set in $B \setminus \{ 0 \}$.
Fix $\e>0$ arbitrarily small. From~\eqref{e:OT}, one infers  $\sup\{ \varphi_k, \, |z| \le 1 - \e \} \le \e$ and $\varphi_k \ge \varphi - \e$
for $k$ large enough. Set $\bar{\varphi}_k = \max \{ 2 \log |z| , \varphi_k - 4 \e\}$ on $\{ |z| \le 1 - \e \}$ and $\bar{\varphi}_k = 2 \log |z|$ on $\{ 1 - \e \le |z| \le 1 \}$. As $ \varphi_k \ge \varphi -\e \ge \log |z|-\e$, $\bar{\varphi}_k$ defines a psh function equal to $\varphi_k - 4 \e$ on the fixed ball $\{ |z|\le e^{-5 \e} \}$. 
Look at $\bar{\varphi}_k + \e \log |z|$: it is a negative psh function on $B$ whose valuative transform equals $\hat{\varphi}_k - \e$. For $k$ sufficiently large,  we get $\hat{\varphi}_k (\nu) - \e \le -1$, and we conclude that $\varphi_k + \e \log |z| \le \varphi \le \varphi_k+ \e$ on $\{ |z|\le e^{-5 \e} \}$ for all $k$ large. This shows uniform convergence on compact subsets of $\varphi_k \to \varphi$ as required.

We have shown that $\varphi$ is a maximal psh weight. Pick a smaller ball $B' \subset B$, and set $\varphi' = \sup\{ u\le 0 \text{ psh on } B' \text{ such that } \hat{u}(\nu) \le -1 \}$. The previous argument shows that  $\varphi$ is a psh weight, $\varphi \ge \log |z| +O(1)$, and
$\hat{\varphi}' = \hat{\varphi}$. 
For suitable constants, 
the psh function $\max\{ \varphi' - C_1, C_2 \log |z|\}$ extends as a negative psh function on $B$, equal to $\varphi' -C_1$ near $0$. It is thus dominated by $\varphi$. The reverse inequality $\varphi \le \varphi'$ is clear by construction. We conclude that the difference $\varphi - \varphi'$ is bounded on $B'$.

Fix a psh germ $u$ defined in a ball $B'$.
Then $u \le -\hat{u}(\nu) \varphi'\le -\hat{u}(\nu) \varphi + O(1)$ by definition, so $\sigma(u,\varphi) \ge -\hat{u}(\varphi)$. On the other hand, $\hat{\varphi}(\nu)= -1$, so $u \le \sigma(u,\varphi) \varphi + O(1)$ implies $-\hat{u}(\nu) \ge \sigma(u,\varphi)$ hence $-\hat{u}(\nu) = \sigma(u,\varphi)$.

 Now if $f\in\cJ(t\varphi)$, then $f\in\elde((t-\e)\hvarphi)$ 
  for each $\e>0$ by Theorem~\ref{thm:mult-comp},
  hence $\log|f|\le t\hvarphi+A$.
  In particular, we get $\nu(f)\ge t-A(\nu)$.
  By definition of $\varphi$,
  this implies that $\log|f|$ is more singular than 
  $(t-C)\varphi$ with $C:=A(\nu)$.
  Hence $\varphi$ is tame.

  We finally prove a strong form of uniqueness: 
  suppose we are given two maximal psh weights
  $\varphi,\varphi'$ satisfying $\sigma(\log|f|,\varphi) =
  \sigma(\log|f|,\varphi')$ for all $f\in \cO$.  Both weights are
  tame and $\hat\varphi= \hat\varphi'$ by what precedes, hence 
  $\varphi =   \varphi'$ up to a bounded function 
  by Theorem~\ref{thm:relat}. This concludes the proof.
\end{proof}

%
%
\subsection{Valuations and generalized Lelong numbers}
Demailly~\cite{dem} defined the \emph{generalized Lelong number} 
$\nu_\varphi(u)$ of
a psh germ $u$ with respect to a psh weight $\varphi$ as the mass at
the origin of the measure $dd^cu\wedge(dd^c\varphi)^{n-1}$. The
generalized Lelong number satisfies the following
\emph{comparison principle}: if $u$
is more singular than $u'$ and $\varphi$ is more singular than
$\varphi'$, then $\nu_{\varphi}(u)\geq\nu_{\varphi'}(u')$. We are now
in position to relate the generalized Lelong number with respect to a tame
weight to an average of valuations, largely 
by repeating the argument given in~\cite{FJ-sing}.
\begin{proof}[Proof of Theorem~B] 
  Pick a tame psh weight $\varphi$ and any psh germ $u$. 
  We first prove that $\hat u$ is integrable with respect to
  the Monge-Amp{\`e}re measure $\MA(\hat\varphi)$.  Let
  $\varphi_k \= \frac1k \Reg(k\varphi)$.  Since $\varphi$ is tame,
  there is a constant $C>0$ such that
  $\varphi\leq\varphi_k\leq(1-C/k)\varphi$
  up to bounded terms. 
  It follows that $Z(\varphi)\leq   Z(\varphi_k)\leq
  (1-C/k)Z(\varphi)$ for each $k$, and thus $\langle
  Z(u),Z(\varphi)^{n-1}\rangle\geq\langle
  Z(u),Z(\varphi_k)^{n-1}\rangle\geq(1-C/k)^{n-1}\langle
  Z(u),Z(\varphi)^{n-1}\rangle$. Since $Z(\varphi_k)$ is a Cartier
  divisor, the Monge-Amp{\`e}re measure $\MA(\hat\varphi_k)$ 
  is finite and $\int_{\cV}\hu\MA(\hat\varphi_k)>-\infty$. 
  It thus follows that $\MA(\hvarphi)$ is also finite.

  The previous argument actually shows 
  $\int_{\cV}\hat u\MA(\hat\varphi)
  =\lim_k \int_{\cV}\hat u\MA(\hat\varphi_k)$.  On the other hand, 
  by the comparison principle, it is clear that
  $\nu_{\varphi}(u)\geq\nu_{\varphi_k}(u)\geq(1-C/k)^{n-1}\nu_{\varphi}(u)$,
  and in particular $\nu_{\varphi_k}(u)$ converges to
  $\nu_{\varphi}(u)$.  We are thus reduced to proving the theorem when
  the weight $\varphi$ has singularities described by a primary ideal
  $\fa$.

  Replacing $\fa$ by an ideal having the same integral closure changes
  the function $\varphi$ by a bounded term only. By choosing generic
  linear combination of a system of generators of $\fa$, we may assume
  that the ideal is generated by exactly $n$ holomorphic 
  functions $f_1,\dots , f_n$.  
  Consider the map $F : (\C^n,0) \to (\C^n,0)$
  given in coordinates by $F = (f_1, \cdots, f_n)$. As $\fa$ 
  is primary, $F$ is finite, and by construction 
  $\varphi = \log | F | $.
  
  Denote by $e(F,\cdot)$ the local topological degree of $F$.
  For any point $p$ near the origin, 
  set $F_*u(p) = \sum_{F(q) = p} e(F,q)u(q)$. 
  Then $F_*u$ is a psh germ, and we have
  \begin{align*}
    \nu_\varphi(u) 
    &=(dd^c u \wedge (dd^c\log |F|)^{n-1})\, \{ 0 \}\\
    &=(dd^c (F_*u) \wedge (dd^c\log |z|)^{n-1})\,  \{ 0 \}
    =\nuL(F_*u).
  \end{align*}

Similarly if $f\in\cO$, then $F_{*}f\in\cO$ is defined by $F_{*}f (p)= \prod
_{F(q)=p} f(q)^{e(F,q)}$. If $\fa$ is primary, then so is $F_{*}\fa$.

  Next, if $u$ has singularities described by an ideal $\fb$, then 
  $\hu=\log|\fb|$, and $F_*u$ has singularities described by $F_*\fb$.
 By the projection formula in intersection theory and~\eqref{e:ideal},
  we get 
  \begin{equation}\label{e:case-log}
    \int_\cV -\hat u\MA(\hat\varphi)
    =\langle \fb, \fa^{\langle n-1 \rangle} \rangle
    =\langle\fb, F^*\fm^{\langle n-1 \rangle} \rangle
    =\langle F_*\fb,\fm^{\langle n-1 \rangle} \rangle
    =\nuL(F_*u),
  \end{equation}
  which shows the formula in Theorem~B in this case.  
  To deal with a general psh germ, 
  we rely again on the approximation technique of Demailly.  
  Pick $\phi \in \cH(u,B)$ with norm $\int_B |\phi|^2 e^{-2u} =1$.
  The mean value theorem yields for all points in the ball $B(r/2)$ of
  radius $r/2$ and centered at the origin:
  \begin{equation*}
    |F_*\phi|^2 
    \le\frac1{c_nr^{2n}}\int_{B(r)} |F_*\phi|^2 \, d\lambda
    \le\frac1{c_nr^{2n}}\, \sup_{B(r)} e^{2F_*u}\cdot
    \int_{B(r)}|F_*\phi|^2e^{-2F_*u}\, d\lambda.
  \end{equation*}
  Now the change of variables formula gives
\begin{align*}
 \int_{B(r)} |F_*\phi|^2 e^{-2F_*u} d\lambda
     &= e(F,0) \cdot \int_{F^{-1} B(r)} \!\!\! |\phi|^2 e^{-2u}|JF|^2 d\lambda \\
 & \le e(F,0)\cdot \sup|JF|^2\cdot\int_{B} |\phi|^2 e^{-2u}d\lambda \le C
\end{align*}
 for some uniform constant $C>0$.  
  Hence $\log |F_*\phi|\le\sup_{B(r)}F_*u-C'\log r$, 
  and by definition of $\Reg(u)$ we infer
  $\sup_{B(r/2)} F_*\Reg(u) \le \sup_{B(r)} F_*u - C' \log r $. 
  Set $u_k\=k^{-1} \Reg{(ku)}$. Then
  \[
  \nuL(F_*u_k) \le \nuL (F_*u) \le \nuL(F_*u_k) + \frac{C''}{k}.
  \]
  In particular we get a sequence of psh functions with singularities
  described by analytical ideals for which 
  $\nuL(F_*u_k)\to\nuL(F_*u)$.
  We saw that $\hat u_k \to \hat u$ pointwise on $\cVqm$. As
  $\MA(\hat\varphi)$ is supported on finitely many divisorial
  valuations (Proposition~\ref{p:suppma}), we infer 
  $\int_\cV-\hu_k\MA(\hat\varphi)\to\int_\cV-\hu\MA(\hat\varphi)$.
  Finally we conclude by applying~\eqref{e:case-log}
  for each $u_k$ and letting $k\to\infty$.
\end{proof}
%
%
\subsection{Proof of Theorem~A}
Consider two psh germs $u,v$. By construction, (1) is equivalent 
to $\hat u=\hat v$. 
Theorem~\ref{thm:repr} shows that (3) implies (1). 
Moreover, (1) implies (3) and (4) by Theorem~\ref{thm:relat} 
and Theorem~B, respectively, and~(1) implies~(2) by
Theorem~\ref{thm:mult-comp}. It remains to be shown that~(2) 
implies~(1). So suppose $\cJ_+ (tu) =\cJ_+(tv)$ for all $t$.
Theorem~\ref{thm:mult-comp} yields 
$\elde (t\hat u)\cap \cO=\elde(t\hat v)\cap \cO$ 
for all $t$, but this does a priori
not imply that $\elde (t \hat u) = \elde(t\hat v)$, 
so we cannot directly apply Proposition~\ref{p:approx} 
to conclude $\hat u=\hat v$. 
However,~\eqref{e:app} applied to $(1+\e)u$ shows that 
$\log|\cJ_+(u)|\ge g_u\ge\log|\cJ_+(u)|-A$, 
and thus $\hat u$ is also the limit of $\frac1k\log|\cJ_+(ku)|$ 
as $k\to\infty$, so the result follows. 

\end{document}